\documentclass[10pt, a4paper]{article}
\usepackage{comment}
\usepackage{graphics}
\usepackage{latexsym}
\usepackage{amssymb}
\usepackage{amsmath}
\usepackage{color}
\usepackage{epsfig}
\usepackage{pgf,tikz}
\usepackage{longtable}
\newcommand{\beq}{\begin{equation}}
\newcommand{\eeq}{\end{equation}}
\newcommand{\beqas}{\begin{eqnarray*}}
\newcommand{\eeqas}{\end{eqnarray*}}
\newcommand{\ep}{\varepsilon}

\newcommand{\ue}{u^{\varepsilon}}

\newcommand{\wc}{\rightharpoonup}

\newcommand{\bY}{{\bf Y}}
\newcommand{\by}{\boldsymbol{y}}

\newcommand{\bV}{{\bf V}}

\newcommand{\bv}{\boldsymbol{v}}
\newcommand{\bu}{\boldsymbol{u}}

\newcommand{\bH}{\boldsymbol{\cal H}}

\newcommand{\fu}{\frak u}
\newcommand{\fv}{\frak v}
\newcommand{\fw}{\frak w}

\newcommand{\bt}{\begin{theorem}}
\newcommand{\bp}{\begin{proposition}}
\newcommand{\bl}{\begin{lemma}}
\newcommand{\br}{\begin{remark}}

\newcommand{\et}{\end{theorem}}
\newcommand{\epr}{\end{proposition}}
\newcommand{\el}{\end{lemma}}
\newcommand{\er}{\end{remark}}

\def\curl{{\rm curl\,}}
\def\scurl{{\rm curl}}

\newcommand{\bproof}{{\it Proof}\ \ }
\newcommand{\eproof}{{\hfill$\Box$}}

\newcommand{\be}{\begin{equation}}
\newcommand{\ee}{\end{equation}}
\newcommand{\bea}{\begin{eqnarray}}
\newcommand{\eea}{\end{eqnarray}}
\newcommand{\beas}{\begin{eqnarray*}}
\newcommand{\eeas}{\end{eqnarray*}}

\def\IR{\mathbb{R}}

\def\curl{{\rm curl\,}}

\def\scurl{{\rm curl}}

\newtheorem{theorem}{Theorem}[section]
\newtheorem{lemma}[theorem]{Lemma}
\newtheorem{assumption}[theorem]{Assumption}
\newtheorem{proposition}[theorem]{Proposition}
\newtheorem{definition}[theorem]{Definition}
\newtheorem{remark}[theorem]{Remark}

\definecolor{darkgreen}{rgb}{0,0.7,0}

\numberwithin{equation}{section}
\numberwithin{equation}{section}
\title{Homogenization of multiscale Maxwell wave equations}
\author{
Van Tiep Chu
and 
Viet Ha Hoang\\[10pt]
Division of Mathematical Sciences,\\
School of Physical and Mathematical Sciences,\\
 Nanyang Technological University, Singapore 637371 
}

\date{}
\begin{document}
\maketitle

\pagestyle{plain}
\begin{abstract}
We study homogenization of multiscale Maxwell wave equation that depends on $n$ separable microscopic scales in a domain $D\subset\IR^d$ on a finite time interval $(0,T)$. Due to the non-compactness of the embedding of $H_0(\curl,D)$ in $L^2(D)^d$, homogenization of Maxwell wave equation can be significantly more complicated than that of scalar wave equations in the $H^1(D)$ setting, and requires analysis uniquely for Maxwell wave equations. We employ multiscale convergence. The homogenized Maxwell wave equation and the initial condition are deduced from the multiscale homogenized equation. When the coefficient of the second order time derivative in the multiscale equation depends on the microscopic scales, the derivation is significantly more complicated, comparing to  scalar wave equations, due to the corrector terms for the solution $\ue$ of the multiscale equation in the $L^2(D)^d$ norm, which do not appear in the $H^1(D)$ setting. For two scale equations, we derive an explicit homogenization error estimate for the case where the solution $u_0$ of the homogenized  equation belongs to $L^\infty((0,T);H^1(\curl,D))$. When $u_0$ only belongs to a weaker regularity space $L^\infty((0,T);H^s(\curl,D))$ for $0<s<1$, we contribute an approach to deduce a new homogenization error in this case, which depends on $s$. For general multiscale problems, a corrector is derived albeit without an explicit homogenization error estimate.  These correctors and homogenization errors play an essential role in deriving numerical correctors for approximating the solutions to the multiscale problems numerically, as considered in our forthcoming publication. 
\end{abstract}

\section{Introduction}
We consider  multiscale Maxwell wave equation \eqref{eq:mseqn} in a domain $D\in\IR^d$ which depends on $n$ microscopic scales. Using multiscale convergence (\cite{Allaire1992}, \cite{Nguetseng1989}, we deduce the multiscale homogenized equation, from which we deduce the homogenized equation. In comparison to multiscale wave equation in the $H^1(D)$ setting, homogenization of a multiscale Maxwell wave equation is quite different due to the embedding of the space $H_0(\curl,D)$ in $L^2(D)^d$ is not compact (see \cite{BLP}). It is far more  complicated to derive the homogenized equation, and the initial condition for the solution of the homogenized equation than for a multiscale scale wave equation when the coefficient $b^\ep$ of the second order time derivative in \eqref{eq:mseqn} depends also on the microscopic scales. This is due to the correctors $\fu_i$ for $\ue$ in $L^2(D)^d$ in \eqref{eq:uemsconv}, and the second order time derivative ${\partial^2\over\partial t^2}\nabla_{y_i}\fu_i$ in \eqref{eq:msprob} which do not belong to $L^2(D\times Y_1\times\ldots\times Y_i)$, and can only be understood in the distribution sense. Derivation of the initial condition requires representations of the time derivative ${\partial\over\partial t}\nabla_{y_i}\fu_i$ in terms of the solutions of the cell problems, which are quite non-trivial. The main contributions of the paper are the new correctors and homogenization errors in two sale problems, especially when the solution to the homogenized equation possesses low regularity, and a general corrector for multiscale problems. Their derivation are quite non trivial, comparing to elliptic and wave equations in the $H^1(D)$ setting.

As is well known, for two scale wave equations, a corrector similar to that of two scale elliptic problems does not hold in general, as the energy of the multiscale problem does not converge to the energy of the homogenized problem (see \cite{Brahim-Otsmane1992}). Therefore, to derive a corrector, we restrict our consideration to the case where the initial condition $g_0$ in \eqref{eq:mseqn} equals zero. For two scale elliptic problems in a domain $D\subset\IR^d$, it is well known that a $O(\ep^{1/2})$ homogenization rate of convergence in the $H^1(D)$ norm can be derived when the solution of the homogenized equation is smooth (see \cite{BLP} and \cite{JKO}). In polygonal domains, which are of interests in numerical discretization, this regularity may not hold. However, if the domain is convex, the solution to the homogenized equation belongs to $H^2(D)$; the $O(\ep^{1/2})$ homogenization rate still holds with the same proof. For wave equation, the $O(\ep^{1/2})$ homogenization rate is deduced in \cite{XHwave} for scalar wave equations and in \cite{XHelasticwave} for elastic wave equations. In this paper, for two scale Maxwell wave equations, we derive the $O(\ep^{1/2})$ homogenization rate when the solution $u_0$ of the homogenized equation belongs to $L^\infty(0,T;H^1(\curl,D))$. However, in polygonal domains, this regularity condition normally does not hold. The solution $u_0$ of the homogenized equation only belongs to a weaker regularity space $L^\infty(0,T;H^s(\curl,D))$ for $0<s<1$. We develop an approach to deduce a new homogenization error for this case of low regularity, using the ideas of \cite{Tiep1}. For multiscale problems, an explicit homogenization error is not available. However, we can deduce a general corrector which requires a procedure quite different from that for multiscale wave equation in \cite{XHwave} due to the corrector terms $\fu_i$ in \eqref{eq:uemsconv}. The correctors and homogenization errors in this paper play a key role in establishing numerical correctors using finite element solutions for the multiscale problems as studied in Chu and Hoang \cite{Tiep3}. 

The paper is organized as follows. In the next section, we define the multiscale Maxwell wave equation, we review the concept of multiscale convergence, extended to functions that depend on the time variable, and use it to derive the multiscale Maxwell wave equation. We then derive the initial condition for the multiscale homogenized Maxwell wave equation, and show that this problem has a unique solution. In Section 3, we derive the homogenized equation from the multicale homogenized equation, together with the initial condition. The derivation is quite non trivial in comparison to the multiscale wave equation in \cite{Brahim-Otsmane1992} and \cite{XHwave} due to the multiscale coefficient $b^\ep$ in \eqref{eq:mseqn} and the corrector terms $\fu_i$ in \eqref{eq:uemsconv}. In Section 4, we study the regularity of the solution $u_0$ of the homogenized equation which is necessary for the derivation of the correctors and homogenization errors. We study correctors in Section 5. For two scale problems, when the solution $u_0$ is in $L^\infty(0,T;H^1(\curl,D))$, we derive the $O(\ep^{1/2})$ homogenization error in the $H(\curl,D)$ norm. When the solution only belongs to a weaker regularity space $L^\infty(0,T;H^s(\curl,D))$ we deduce a weaker homogenization convergence rate. The proof requires substantial technical developments. For general multiscale problems, we derive in this section a general corrector, without an explicit convergence rate. 

Throughout the paper, without indicating a variable, $\curl$ and $\nabla$  denote the $\curl$ and gradient of a function of the variable $x$, with respect to $x$, and by $\curl_x$ and $\nabla_x$, we denote partial $\curl$ and partial gradient with respect to $x$, of a function depending on $x$ and other variables. Repeated indices indicate summation.

\section{Multiscale homogenization of multiscale Maxwell wave equation}
\subsection{Problem setting}
Let $D$ be a bounded domain in $\IR^d$ ($d=2,3$). Let $Y$ be the unit cube in $\IR^d$. By $Y_1,\ldots,Y_n$ we denote $n$ copies of $Y$. We denote by $\bY$ the product set $Y_1\times Y_2\times\ldots\times Y_n$ and by $\by=(y_1,\ldots,y_n)$. For $i=1,\ldots,n$, we denote by $\bY_i=Y_1\times\ldots\times Y_i$. Let $a$ and $b$ be  functions from $D\times Y_1\times\ldots\times Y_n$ to $\IR^{d\times d}_{sym}$. We assume that the symmetric matrix functions $a$ and $b$ satisfy the boundedness and coerciveness conditions: for all $x\in D$ and $\by\in\bY$, and all $\xi,\zeta\in \IR^d$,
\be
\begin{array}{lr}
\displaystyle{\alpha|\xi|^2\le a_{ij}(x,\by)\xi_i\xi_j,\ \ \ a_{ij}\xi_i\zeta_j\le\beta|\xi||\zeta|}\\
\displaystyle{\alpha|\xi|^2\le b_{ij}(x,\by)\xi_i\xi_j,\ \ \ b_{ij}\xi_i\zeta_j\le\beta|\xi||\zeta|}
\end{array}
\label{eq:coercive}
\ee
where $\alpha$ and $\beta$ are positive numbers. 
 Let $\ep$ be a small positive value. Let $\ep_1,\ldots,\ep_n$ be $n$ functions of $\ep$ that denote the $n$ microscopic scales that the problem depends on. We assume the following scale separation properties: for all $i=1,\ldots,n-1$
\be
\lim_{\ep\to 0}{\ep_{i+1}(\ep)\over\ep_i(\ep)}=0.
\label{eq:scaleseparation}
\ee
Without loss of generality, we assume that $\ep_1(\ep)=\ep$. 
We define the multiscale coefficient of the Maxwell equation $a^\ep$ and $b^\ep$ which are functions from $D$ to $\IR^{d\times d}_{sym}$ as 
\[
a^\ep(x)=a(x,{x\over\ep_1},\ldots,{x\over\ep_n}),\ \ b^\ep(x)=b(x,{x\over\ep_1},\ldots,{x\over\ep_n}).
\]
When $d=3$ we define the spaces
\[
W=H_0(\curl,D)=\{u\in L^2(\Omega)^3,\ \ \curl u\in L^2(\Omega)^3,\ \ u\times\nu=0\},\ \ 
H=L^2(D)^3
\]
and when $d=2$
\[
W=H_0(\curl,D)=\{u\in L^2(D)^2,\ \ \curl u\in L^2(D),\ \ u\times\nu=0\},\ \ 
H=L^2(D)^2
\]
where $\nu$ denotes the outward normal vector on the boundary $\partial D$. We have the Gelfand triple $W\subset H\subset W'$. We denote by $\langle\cdot,\cdot\rangle$ the inner product in $H$, extending to the duality pairing between $W'$ and $W$. 
 We note that when $d=3$, $\curl\ue$ is a vector function in $L^2(D)^3$ and when $d=2$, $\curl\ue$ is a scalar function in $L^2(D)$. 
Let $f\in L^2(0,T; H)$, $g^0\in W$ and $g^1\in H$. 
We consider the problem:
Find $u^\varepsilon (t,x) \in L^2(0,T; W) $
\be
\begin{cases}
b^\varepsilon(x)\dfrac{\partial^2u^\varepsilon(t,x) }{\partial t^2}+ \mathrm{curl}( a^\varepsilon(x)\mathrm{curl} u^\varepsilon(t,x) )=f(t,x),&  (t,x)\in (0,T) \times D\\
u^\varepsilon(0,x)=g_0(x)&\\
u_t^\varepsilon(0,x)=g_1(x)&\\
\end{cases}
\label{eq:mseqn}
\ee
with the boundary condition $\ue\times\nu=0$ on $\partial D$.  We will mostly present the analysis for the case $d=3$ and only discuss the case $d=2$ when there is significant difference.  For notational conciseness, we denote by
\be
H=L^2(D)^3,\ \ H_i=L^2(D\times \bY_i)^3,\ \ i=1,\ldots,n.
\label{eq:H}
\ee
In variational form, this problem becomes: Find $\ue\in L^2(0,T; W) \cap H^1(0,T;H)$ so that
\be
\left\langle b^\varepsilon(x)\frac{\partial^2 \ue}{\partial t^2}, \phi(x)\right\rangle_{W',W}+\int_Da^\varepsilon(x)\mathrm{curl} u^\varepsilon(t,x)\cdot\curl\phi(x)dx
=\int_Df(t,x)\cdot\phi(x)dx
\label{eq:varprob}
\ee
for all $\phi\in W$ when $d=3$; and when $d=2$ we need to replace the vector product for $\curl$ by the scalar multiplication. 
Problem \eqref{eq:varprob} has a unique solution $\ue\in L^2(0,T;W)\bigcap H^1(0,T;H)\bigcap H^2(0,T;W')$ that satisfies
\be
\|\ue\|_{L^2(0,T;W)}+\|\ue\|_{H^1(0,T;H)}+\|\ue\|_{H^2(0,T;W')}\le c(\|f\|_{L^2(0,T;H)}+\|g^0\|_W+\|g^1\|_H)
\label{eq:uepbound}
\ee
where the constant $c$ only depends on the constants $\alpha$ and $\beta$ in \eqref{eq:coercive} and $T$ (see Wloka \cite{Wloka}). 

We will study this problem via multiscale convergence.

\subsection{Multiscale convergence}
We recall the definition of multiscale convergence (see Nguetseng \cite{Nguetseng1989}, Allaire \cite{Allaire1992} and Allaire and Briane \cite{Allaire1996}. 
\begin{definition}
A sequence of functions $\{w^\ep\}_\ep\subset L^2(0,T;H)$ $(n+1)$-scale converges to a function $w_0\in L^2((0,T);D\times \bY)$ if for all smooth functions $\phi(t,x, \by)$ which are $Y$ periodic w.r.t. $y_i$ for all $i=1,\ldots,n$:
\[
\lim_{\ep \to 0}\int_0^T\int_Dw^\ep(t,x)\phi(t,x,{x\over\ep_1},\ldots,{x\over\ep_n})dxdt=\int_0^T\int_D\int_\bY w_0(t,x,\by)\phi(t,x,\by)d\by dxdt.
\]
\end{definition}
We have the following result.
\begin{proposition}\label{prop:2scexistence}
From a bounded sequence in $L^2(0,T; H)$ we can extract an $(n+1)$-scale convergent subsequence.
\end{proposition}
We note that the definition above for functions which depend also on $t$ is slightly different from that in \cite{Nguetseng1989} and \cite{Allaire1992} as we take also the integral with respect to $t$. However, the proof of proposition \ref{prop:2scexistence} is similar.

For a bounded sequence in $L^2(0,T;W)$, we have the following results  which are very similar to those in \cite{Tiep1} for functions which do not depend on $t$. The proofs for these results follow the same lines of those in \cite{Tiep1} so we do not present them here. As in \cite{Wellander} and \cite{Tiep1}, we denote by $\tilde H_\#(\curl,Y_i)$ the space of equivalent classes of functions in $H_\#(\curl,Y_i)$ of equal $\curl$. 
\begin{proposition}\label{prop:msHcurl}
Let $\{w^\ep\}_\ep$ be a bounded sequence in  $L^2(0,T; W)$. There is a subsequence (not renumbered), a function $w_0\in  L^2(0,T; W)$, $n$ functions $\frak w_i\in L^2((0,T) \times D\times Y_1\times\ldots\times Y_{i-1},H^1_\#(Y_i)/\IR)$ such that 
\[
w^\ep \stackrel{(n+1)-\mbox{scale}}{\longrightarrow}w_0+\sum_{i=1}^n\nabla_{y_i}\frak w_i.
\]
Further, there are $n$ functions $w_i\in L^2((0,T) \times D\times\ldots\times Y_{i-1},\tilde H_\#(\curl,Y_i))$ such that
\[
\curl w^\ep \stackrel{(n+1)-\mbox{scale}}{\longrightarrow}\curl w_0+\sum_{i=1}^n\scurl_{y_i}w_i.
\]
\end{proposition}

\subsection{Multiscale homogenized Maxwell wave problem}
From \eqref{eq:uepbound} and Proposition \eqref{prop:msHcurl}, we can extract a subsequence (not renumbered), a function $u_0\in L^2(0,T; W)$, $n$ functions $\frak u_i\in L^2(0,T;D\times Y_1\times\ldots\times Y_{i-1},H^1_\#(Y_i)/\IR)$ and $n$ functions $u_i\in L^2(0,T;D\times Y_1\times\ldots\times Y_{i-1},\tilde H_\#(\curl,Y_i))$ such that
\be
\ue\stackrel{(n+1)-\mbox{scale}}{\longrightarrow} u_0+\sum_{i=1}^n\nabla_{y_i}\frak u_i,
\label{eq:uemsconv}
\ee
and 
\be
\curl\ue \stackrel{(n+1)-\mbox{scale}}{\longrightarrow} \curl u_0+\sum_{i=1}^n\scurl_{y_i}u_i.
\label{eq:curluemsconv}
\ee
For $i=1,\ldots,n$, let $W_i=L^2(D\times Y_1\times\ldots\times Y_{i-1},\tilde H_\#(\curl,Y_i))$ and $V_i=L^2(D\times Y_1\times\ldots\times Y_{i-1},H^1_\#(Y_i)/\IR)$. 
We define the space $\bV$
\[
\bV=W\times W_1\times\ldots\times W_n\times V_1\times\ldots\times V_n.
\]
For $\bv=(v_0,\{v_i\},\{\frak v_i\})\in \bV$, we define the form
\[
|||\bv|||=\|v_0\|_{H(\curl,D)}+\sum_{i=1}^n\|v_i\|_{L^2(D\times\bY_{i-1},\tilde H_\#(\curl,Y_i))}+\sum_{i=1}^n\|\frak v_i\|_{L^2(D\times\bY_{i-1},H^1_\#(Y_i))}.
\]
Let $\bu=(u_0,\{u_i\}, \{\frak u_i\}) \in \bV$. We define the function
\[
\int_{\bY}b(x,\by)\left(u_0(t,x)+\sum_{i=1}^n\nabla_{y_i}\fu_i(t,x,\by)\right)d\by
\]
in $W'$ as
\[
\left\langle\int_{\bY}b(x,\by)\left(u_0(t,x)+\sum_{i=1}^n\nabla_{y_i}\fu_i(t,x,\by)\right)d\by,v_0\right\rangle_{W',W}=\int_D\int_\bY b(x,\by)\left(u_0(t,x)+\sum_{i=1}^n\nabla_{y_i}\fu_i(t,x,\by)\right)\cdot v_0d\by dx;
\]
and the function $b(x,\by)(u_0(t,x)+\sum_{i=1}^n\nabla_{y_i}\fu_i(t,x,\by))$ in $V_j'$ as
\[
\left\langle b(x,\by)\left(u_0(t,x)+\sum_{i=1}^n\nabla_{y_i}\fu_i(t,x,\by)\right),\fv_j\right\rangle_{V_j',V_j}=\int_D\int_\bY b(x,\by)(u_0(t,x)+\sum_{i=1}^n\nabla_{y_i}\fu_i(t,x,\by))\cdot\nabla_{y_j}\fv_j(x,\by_j) d\by dx.
\]
We then have the following result.
\begin{proposition}\label{prop:2.4}
The function $\bu=(u_0,\{u_i\},\{\bu_i\})$ satisfies
\beqas
&&\left\langle{\partial^2\over\partial t^2}\int_{\bY}b(x,\by)\left(u_0(t,x)+\sum_{i=1}^n\nabla_{y_i}\fu_i(t,x,\by)\right)d\by,v_0\right\rangle_{W',W}=\\
&&\int_D f(x)\cdot v_0(x)dx-\int_D\int_\bY a(x,\by)\left(\curl u_0+\sum_{i=1}^n\scurl_{y_i}u_i\right)\cdot\left( \curl v_0+\sum_{i=1}^n\scurl_{y_i}v_i\right)d\by dx
\eeqas
and 
\[
\left\langle{\partial^2\over\partial t^2}\left(b(x,\by)\left(u_0(t,x)+\sum_{i=1}^n\nabla_{y_i}\fu_i(t,x,\by)\right)\right),\fv_j\right\rangle_{W_j',W_j}=0,\ \ j=1,\ldots,n,
\]
i.e. the function $\bu$ satisfies the multiscale homogenized equation
\begin{multline}
\left\langle{\partial^2\over\partial t^2}\int_\bY b(x,\by)\left(u_0(t,x)+\sum_{i=1}^n\nabla_{y_i}\fu_i\right)d\by,v_0\right\rangle_{W',W}+\\
\sum_{j=1}^n\left({\partial^2\over\partial t^2}b(x,\by)\left(u_0(t,x)+\sum_{i=1}^n\nabla_{y_i} \frak u_i(t,x,\by)\right),\frak v_j\right)_{V_j',V_j}\\
+\int_D\int_\bY a(x,\by)\left(\curl u_0+\sum_{i=1}^n\scurl_{y_i}u_i\right)\cdot\left( \curl v_0+\sum_{i=1}^n\scurl_{y_i}v_i\right)d\by dx=\int_D f(x)\cdot v_0(x)dx
\label{eq:msprob}
\end{multline}
for all $\bv=(v_0,\{v_i\},\{\frak v_i\})\in \bV$. 
\end{proposition}
\bproof
Let $q\in {\cal D}(0,T)$. Let $v_0\in {\cal D}(D)$, $v_i\in {\cal D}(D,C^\infty_\#(Y_1,\ldots,C^\infty_\#(Y_i)\ldots))^3$ and $\fv_i\in {\cal D}(D,C^\infty_\#(Y_1,\ldots,C^\infty_\#(Y_i)\ldots))$ for $i=1,\ldots,n$. Choosing a test function of the form
$$
\phi(t,x)= \left(v_0(x)+ \sum_{i=1}^n\ep_i v_i(x,\frac{x}{\ep_1},\ldots,\frac{x}{\ep_i})+\sum_{i=1}^n\ep_i \nabla{\fv}_i(x,\frac{x}{\ep_1},\ldots,\frac{x}{\ep_i})\right) q(t)
$$
we obtain
\beqas
&&\int_0^T \int_D b(x,\frac{x}{\ep_1},\ldots,\frac{x}{\ep_n})u^\ep(t,x) \cdot \left(v_0(x)+ \sum_{i=1}^n\ep_i v_i(x,\frac{x}{\ep_1},\ldots,\frac{x}{\ep_i})+\right.\\
&&\qquad\qquad\left.\sum_{i=1}^n\left(\ep_i\nabla_x\fv_i(x,\frac{x}{\ep_1},\ldots,\frac{x}{\ep_i})+\sum_{j=1}^i{\ep_i\over\ep_j} \nabla_{y_j}\fv_i(x,\frac{x}{\ep_1},\ldots,\frac{x}{\ep_i})\right)\right) q''(t) dxdt\\
&&+\int_0^T \int\limits_Da(x,{x\over\ep_1},\ldots,{x\over\ep_n}) \curl\ue(t,x)\cdot \\
&&\qquad\qquad\left(\curl v_0(x)+\sum_{i=1}^n\left(\ep_i\scurl_x v_i(x,\frac{x}{\ep_1},\ldots,\frac{x}{\ep_i})+\sum_{j=1}^i{\ep_i\over\ep_j}\scurl_{y_j} v_i(x,\frac{x}{\ep_1},\ldots,\frac{x}{\ep_i})\right)\right) q(t)dxdt\\
&&=\int_0^T \int_Df(t,x)\cdot \left(v_0(x)+ \sum_{i=1}^n\ep_i v_i(x,\frac{x}{\ep_1},\ldots,\frac{x}{\ep_i})+\sum_{i=1}^n\ep_i \nabla_x\fv_i(x,\frac{x}{\ep_1},\ldots,\frac{x}{\ep_i})+\right.\\
&&\qquad\qquad\left.\sum_{i=1}^n\sum_{j=1}^i{\ep_i\over\ep_j}\nabla_{y_j}\fv_i(x,\frac{x}{\ep_1},\ldots,\frac{x}{\ep_i})\right) q(t)dxdt
\eeqas
Passing to the two scale limit, using the scale separation \eqref{eq:scaleseparation}, we have
\beqas
&&\int_0^T \int_D\int_\bY b(x,\by)\left(u_0(t,x)+\sum_{i=1}^n\nabla_{y_i}\fu_i(t,x,\by_i)\right) \cdot \left(v_0(x)+ \sum_{i=1}^n\nabla_{y_i}\fv_i(x,\by_i)\right) q''(t)d\by dxdt\\
&&+\int_0^T \int_D\int_{\bY}a(x,\by) \left(\curl u_0(t,x)+\sum_{i=1}^n\curl_{y_i} u_i(t,x,y)\right)\cdot\left(\curl v_0(x) +\sum_{i=1}^n\scurl_{y_i} v_i(x,\by_i)\right) q(t)d\by dxdt\\
&&=\int_0^T \int_D\int_\bY f(t,x)\cdot \left(v_0(x)+\sum_{i=1}^n\nabla_{y_i}\fv_i(x,\by_i)\right) q(t)d\by dxdt\\
&&=\int_0^T\int_Df(t,x)\cdot v_0(x)q(t)dxdt.
\eeqas
Using a density argument, we find that this equation holds for all $(v_0,\{v_i\},\{\fv_i\})\in\bV$. The conclusion then follows. 
\eproof

We now establish the initial conditions.
\bp We have $u_0\in H^1(0,T;H)$, $\nabla_{y_i}\frak u_i\in H^1(0,T;H_i)$ for all $i=1,\ldots,n$. Further
\be
u_0(0,\cdot)=g_0,\ \ \nabla_{y_i}\frak u_i(0,\cdot,\cdot)=0,
\label{eq:u0i}
\ee
\be
{\partial\over\partial t}\int_\bY b(x,\by)(u_0(t,x)+\sum_{j=1}^n\nabla_{y_j}\fu_j(t,x,\by_j))d\by\bigg|_{t=0}=\int_\bY b(x,\by)g_1(x)d\by,\ \ \mbox{in}\ W'
\label{eq:fuiv0i}
\ee
and for $i=1,\ldots,n$
\be
{\partial\over\partial t}b(x,\by)\left(u_0(t,x)+\sum_{j=1}^n\nabla_{y_j}\fu_j(t,x,\by)\right)\bigg|_{t=0}=b(x,\by)g_1(x),\ \ \mbox{in}\ V_i'.
\label{eq:fuivii}
\ee
\epr
\bproof
As $\ue$ is bounded in $H^1(0,T;H)$, $u_0$ belongs to $H^1(0,T;H)\subset C([0,T];H)$ and is the weak limit in this space of $\ue$. Let $\phi\in C^\infty([0,T]\times D)$ with $\phi=0$ when $t=T$. We have
\beqas
\int_0^T\int_D{\partial\ue\over\partial t}\phi(t,x)dxdt=\int_0^T\int_D\left({\partial\over\partial t}(\ue(t,x)\phi(t,x))-\ue(t,x){\partial\phi\over\partial t}(t,x)\right)dxdt\\
=-\int_D\ue(0,x)\phi(0,x)-\int_0^T\int_D\ue(t,x){\partial\phi\over\partial t}(t,x)dxdt\\
\to -\int_Dg_0(x)\phi(0,x)dx-\int_0^T\int_Du_0(t,x){\partial\phi\over\partial t}(t,x)dxdt.
\eeqas
On the other hand
\beqas
\int_0^T\int_D{\partial\ue\over\partial t}\phi(t,x)dxdt\to\int_0^T\int_D{\partial u_0\over\partial t}(t,x)\phi(t,x)=-\int_Du_0(0,x)\phi(0,x)-\int_0^T\int_Du_0(t,x){\partial\phi\over\partial t}(t,x)dxdt.
\eeqas
Thus $u_0(0,x)=g_0$. 

As $\{{\partial\ue\over\partial t}\}_\ep$ is bounded in $L^2(0,T;H)$ so there is a subsequence that $(n+1)$-scale converges. Let $\xi\in L^2((0,T);L^2(D\times\bY))$ be the $(n+1)$-scale limit. Let $\phi(t,x,y_1)\in C^\infty_0((0,T);C^\infty_0 (D,C^\infty_\#(Y_1)))$. We have that
\[
\lim_{n\to\infty}\int_0^T\int_D{\partial\ue\over\partial t}(t,x)\phi(t,x,{x\over\ep_1})dxdt\to\int_0^T\int_D\int_{\bY}\xi(t,x,\by)\phi(t,x,y_1)d\by dxdt.
\]
On the other hand 
\beqas
\int_0^T\int_D{\partial\ue\over\partial t}(t,x)\phi(t,x,{x\over\ep_1})dxdt=-\int_0^T\int_D\ue(t,x){\partial\phi\over\partial t}(t,x,{x\over\ep_1})dxdt
\eeqas
which converges to 
\[
-\int_0^T\int_D\int_{\bY}(u_0+\sum_{i=1}^n\nabla_{y_i}\fu_i){\partial\phi\over\partial t}(t,x,y_1)d\by dxdt=-\int_0^T\int_D\int_{Y_1}(u_0+\nabla_{y_1}\fu_1){\partial\phi\over\partial t}(t,x,y_1)dy_1dxdt.
\]
Thus 
\[
\int_{Y_2}\ldots\int_{Y_n}\xi(t,x,\by)dy_n\ldots dy_2={\partial\over\partial t}(u_0+\nabla_{y_1}\fu_1)
\]
so 
\[
{\partial\over\partial t}\nabla_{y_1}\fu_1=\int_{Y_2}\ldots\int_{Y_n}\xi(t,x,\by)dy_n\ldots dy_2-{\partial u_0\over\partial t}\in H_1
\]
(we refer to \eqref{eq:H} for the definition of the spaces $H_i$). 
Similarly, using a function $\phi(t,x,y_1,y_2)\in C^\infty_0((0,T),C^\infty_0 (D,C^\infty_\#(Y_1,C^\infty_\#(Y_2))))$, we have
\[
{\partial\over\partial t}\nabla_{y_2}\fu_2=\int_{Y_3}\ldots\int_{Y_n}\xi(t,x,\by)dy_n\ldots dy_3-{\partial u_0\over\partial t}-{\partial\over\partial t}\nabla_{y_1}\fu_1\in H_2.
\]
Continuing this process, we have that for all $i=1,\ldots,n$,
\[
{\partial\over\partial t}\nabla_{y_i}\fu_i\in H_i
\]
Finally, we have
\[
\xi(t,x,\by)={\partial u_0\over\partial t}+\sum_{i=1}^n{\partial\over\partial t}\nabla_{y_i}\fu_i.
\]
Let $q\in C^\infty([0,T])$ with $q(T)=0$. Let $\phi(x)=v_0(x)+\sum_{i=1}^n\ep_iv_i(x,{x\over\ep_1},\ldots,{x\over\ep_i})+\sum_{i=1}^n\ep_i\nabla\fv_i(x,{x\over\ep_1},\ldots,{x\over\ep_i})$. We have
\begin{eqnarray}
&&\int_0^T\left\langle b^\ep(x){\partial^2\ue\over\partial t^2},\phi\right\rangle_{W',W} q(t)dt=\int_0^T{\partial\over\partial t}\left(\left\langle b^\ep{\partial\ue\over\partial t},\phi\right\rangle_{W',W} q(t)\right)dt -\int_0^T\left\langle b^\ep{\partial\ue\over\partial t},\phi\right\rangle_{W',W} {dq(t)\over dt}dt\nonumber\\
&&=-\left\langle b^\ep{\partial\ue\over\partial t}(0),\phi\right\rangle_{W',W} q(0)-\int_0^T\left\langle b^\ep{\partial\ue\over\partial t},\phi\right\rangle_{W',W}{dq(t)\over dt}dt\nonumber\\
&&=-\langle b^\ep g^1,\phi\rangle q(0)-\int_0^T\langle b^\ep{\partial \ue\over\partial t},\phi\rangle_{W',W} {dq(t)\over dt}dt\nonumber\\
&&\to -\int_D\int_\bY b(x,\by)g_1(x)(v_0(x)+\sum_{i=1}^n\nabla_{y_i}\fv_i(x,\by_i))q(0)d\by dx-\nonumber\\
&&\int_0^T\int_D\int_Yb(x,\by)\left({\partial u_0\over\partial t}(x)+\sum_{i=1}^n{\partial\over\partial t}\nabla_{y_i}\fu_i(t,x,\by_i)\right)\cdot\left(v_0(x)+\sum_{i=1}^n\nabla_{y_i}\fv_i(x,\by_i)\right){dq(t)\over dt}d\by dxdt
\label{eq:sec2i1}
\end{eqnarray}
when $\ep\to 0$. 
On the other hand, let $q_n$ be a sequence in $C^\infty_0(0,T)$ that converges to $q(t)$ in $L^2(0,T)$ when $n\to \infty$. As $b^\ep{\partial^2\over\partial t^2}\ue$ is bounded in $L^2(0,T;W')$ so  there is a constant $c>0$ such that
\beqas
\left|\int_0^T\langle b^\ep{\partial^2\ue\over\partial t^2},\phi\rangle_{W',W} q_n(t)dt-\int_0^T\langle b^\ep{\partial^2\ue\over\partial t^2},\phi\rangle_{W',W} q(t)dt\right|\le c\|q_n-q\|_{L^2(0,T)}.
\eeqas 
As $q_n\in C^\infty_0(0,T)$, when $\ep\to0$, 
\beqas
&&\lim_{\ep\to 0}\int_0^T\langle b^\ep{\partial^2\ue\over\partial t^2},\phi\rangle q_n(t)dt=\\
&&\int_0^T\int_D\int_{\bY}b(x,\by)(u_0(t,x)+\sum_{i=1}^n\nabla_{y_i}\fu_i(t,x,\by_i))\cdot(v_0(x)+\sum_{i=1}^n\nabla_{y_i}\fu_i(x,\by_i)\rangle q_n''(t)d\by dxdt\\
&&=\int_0^T\langle{\partial^2\over\partial t^2}\int_{\bY}b(x,\by)(u_0(t,x)+\sum_{i=1}^n\nabla_{y_i}\fu_i(t,x,\by_i))d\by,v_0\rangle_{W',W}q_n(t)dt+\\
&&\sum_{i=1}^n\int_0^T\langle{\partial^2\over\partial t^2}b(x,y)(u_0(t,x)+\sum_{j=1}^n\nabla_{y_j}\fu_j(t,x,\by_i),\fv_i\rangle_{V_i',V_i}q_n(t)dt.
\eeqas
Passing to the limit when $n\to\infty$, we have
\begin{eqnarray}
&&\lim_{\ep\to 0}\int_0^T\langle b^\ep{\partial^2\ue\over\partial t^2},\phi\rangle q(t)dt=\nonumber\\
&&\int_0^T\langle{\partial^2\over\partial t^2}\int_{\bY}b(x,\by)(u_0(t,x)+\sum_{i=1}^n\nabla_{y_i}\fu_i(t,x,\by_i))d\by,v_0\rangle_{W',W}q(t)dt+\nonumber\\
&&\qquad\qquad\sum_{i=1}^n\int_0^T\langle{\partial^2\over\partial t^2}b(x,y)(u_0(t,x)+\sum_{j=1}^n\nabla_{y_j}\fu_j(t,x,\by_j),\fv_i\rangle_{V_i',V_i}q(t)dt.
\label{eq:expr}
\end{eqnarray}
From Proposition \ref{prop:2.4}, as ${\partial^2\over\partial t^2}\int_{\bY}b(x,\by)(u_0(t,x)+\sum_{j=1}^n\nabla_{y_j}\fu_j(t,x,\by_i))d\by\in L^2(0,T;W')$, ${\partial\over\partial t}\int_{\bY}b(x,\by)(u_0(t,x)+\sum_{i=1}^n\nabla_{y_i}\fu_i(t,x,\by_i))d\by\in C([0,T];W')$ so the initial condition ${\partial\over\partial t}\int_{\bY}b(x,\by)(u_0(t,x)+\sum_{j=1}^n\nabla_{y_j}\fu_j(t,x,\by_j))d\by$ at $t=0$ is well defined in $W'$. Similarly, the initial condition ${\partial\over\partial t}b(x,\by)(u_0+\sum_{j=1}^n\nabla_{y_j}\fu_i(t,x,\by_j)$ at $t=0$ is well defined in $V_i'$. 
The right hand side of \eqref{eq:expr} can be written as
\begin{eqnarray}
&&\int_0^T{\partial\over\partial t}\left(\left\langle{\partial\over\partial t}\int_{\bY}b(x,\by)(u_0(t,x)+\sum_{j=1}^n\nabla_{y_j}\fu_j(t,x,\by_i))d\by,v_0\right\rangle_{W',W}q(t)\right)dt-\nonumber\\
&&\int_0^T\left\langle{\partial\over\partial t}\int_{\bY}b(x,\by)(u_0(t,x)+\sum_{j=1}^n\nabla_{y_j}\fu_j(t,x,\by_j)dy,v_0\right\rangle_{W',W}{dq(t)\over dt}dt+\nonumber\\
&&\sum_{i=1}^n\int_0^T{\partial\over\partial t}\left(\left\langle{\partial\over\partial t}(b(x,\by)(u_0(t,x)+\sum_{j=1}^n\nabla_{y_j}\fu_j(t,x,\by_j)),\fv_i\right\rangle_{V_j',V_j}q(t)\right)dt-\nonumber\\
&&\sum_{i=1}^n\int_0^T\left\langle{\partial\over\partial t}(b(x,\by)(u_0(t,x)+\sum_{j=1}^n\nabla_{y_j}\fu_j(t,x,\by_j)),\fv_i\right\rangle_{V_i',V_i}{dq(t)\over dt}dt\nonumber\\
&&=-\left\langle{\partial\over\partial t}\int_\bY b(x,\by)(u_0+\sum_{j=1}^n\nabla_{y_j}\fu_j)d\by\bigg|_{t=0},v_0\right\rangle_{W',W}q(0)-\nonumber\\
&&\int_0^T\left\langle{\partial\over\partial t}\int_{\bY}b(x,\by)(u_0(t,x)+\sum_{j=1}^n\nabla_{y_j}\fu_j(t,x,\by_j)dy,v_0\right\rangle_{W',W}{dq(t)\over dt}dt-\nonumber\\
&&\sum_{i=1}^n\left\langle {\partial\over\partial t}\left(b(x,\by)(u_0+\sum_{j=1}^n\nabla_{y_j}\fu_j)\bigg|_{t=0}\right),\fv_i\right\rangle_{V_i',V_i}q(0)\nonumber\\
&&-\sum_{i=1}^n\int_0^T\left\langle{\partial\over\partial t}(b(x,\by)(u_0(t,x)+\sum_{j=1}^n\nabla_{y_j}\fu_j(t,x,\by_j)),\fv_i\right\rangle_{V_i',V_i}{dq(t)\over dt}dt.
\label{eq:sec2i2}
\end{eqnarray}
Comparing \eqref{eq:sec2i1} and \eqref{eq:sec2i2}, we have
\[
\left\langle{\partial\over\partial t}\int_\bY b(x,\by)(u_0(t,x)+\sum_{i=1}^n\nabla_{y_j}\fu_j(t,x,\by_j)d\by\Big|_{t=0},v_0\right\rangle_{W',W}=\int_D\int_\bY b(x,\by)g_1(x)\cdot v_0(x)d\by dx,
\]
and 
\[
\left\langle{\partial\over\partial t}(b(x,\by)(u_0+\sum_{j=1}^n\nabla_{y_j}\fu_j)\bigg|_{t=0},\fv_i\right\rangle_{V_i',V_i}=\int_D\int_\bY b(x,\by)g_1(x)\cdot\nabla_{y_i}\fv_i(x,\by_i)d\by dx.
\]
for all $i=1,\ldots,n$. We get the desired result.

\eproof

\bp With the initial conditions \eqref{eq:u0i}, \eqref{eq:fuiv0i} and \eqref{eq:fuivii}, problem \eqref{eq:msprob} has a unique solution.
\epr
\bproof
We show that when $f=0$, $g_0=0$ and $g_1=0$, the solution of \eqref{eq:msprob} is $u_0=0$, $u_i=0$ and $\fu_i=0$ for all $i=1,\ldots,n$.

Following the procedure in \cite{Wloka} Theorem 19.1 for showing the uniqueness of a solution of a wave equation, fixing $s \in (0,T)$, we define 
\begin{eqnarray}
  w_0(t)=
  \begin{cases}
\displaystyle  -\int_t^s u_0(\sigma) d\sigma, &t < s,\\
  0,& t\geq s;
  \end{cases}\nonumber\\
 \frak w_1(t)=
  \begin{cases}
\displaystyle  -\int_t^s \frak u_1(\sigma) d\sigma, &t < s,\\
  0,& t\geq s;
  \end{cases}\label{eq:w}\\
  w_1(t)=
  \begin{cases}
\displaystyle  -\int_t^s u_1(\sigma) d\sigma, &t < s,\nonumber\\
  0,& t\geq s.
  \end{cases}
 \end{eqnarray}
We have
\beqas
&&{\partial\over\partial t}\left\langle{\partial\over\partial t}\int_\bY b(x,\by)(u_0(t,x)+\sum_{i=1}^n\nabla_{y_i}\fu_i(t,x,\by_i)d\by,w_0(t,\cdot)\right\rangle_{W',W}\\
&&\qquad\qquad+{\partial\over\partial t}\left\langle{\partial\over\partial t}b(x,\by)(u_0(t,x)+\sum_{i=1}^n\nabla_{y_i}\fu_i(t,x,\by_i)),\fw_i(t,\cdot,\cdot)\right\rangle_{V_i',V_i}\\
&&\qquad\qquad=\left\langle {\partial^2\over\partial t^2}\int_\bY b(x,\by)(u_0+\sum_{i=1}^n\nabla_{y_i}\fu_i)d\by,w_0\right\rangle_{W',W}
+\left\langle{\partial^2\over\partial t^2}b(x,\by)(u_0+\sum_{i=1}^n\nabla_{y_i}\fu_i),\fw_i\right\rangle_{V_i',V_i}\\
&&\qquad\qquad+\int_D\int_\bY b(x,\by)\left({\partial u_0\over\partial t}+\sum_{i=1}^n{\partial\over\partial t}\nabla_{y_i}\fu_i\right)\cdot\left({\partial w_0\over\partial t}+\sum_{i=1}^n{\partial\over\partial t}\nabla_{y_i}\fw_i\right)d\by dx\\
&&\qquad\qquad=-\int_D\int_\bY a(x,\by)\left(\curl u_0(t,x)+\sum_{i=1}^n\scurl_{y_i}u_i(t,x,y_i)\right)\cdot\left(\curl w_0(t,x)+\sum_{i=1}^n\scurl_{y_i}w_i(t,x,y)\right)d\by dx\\
&&\qquad\qquad+\int_D\int_\bY b(x,\by)\left({\partial u_0\over\partial t}+\sum_{i=1}^n{\partial\over\partial t}\nabla_{y_i}\fu_i\right)\cdot\left({\partial w_0\over\partial t}+\sum_{i=1}^n{\partial\over\partial t}\nabla_{y_i}\fw_i\right)d\by dx.
\eeqas
Integrating over $(0,s)$, we get
\beqas
&&\left\langle{\partial\over\partial t}\int_\bY b(x,\by)(u_0(t,x)+\sum_{i=1}^n\nabla_{y_i}\fu_i(t,x,\by_i))d\by\Big|_{t=s},w_0(s,\cdot)\right\rangle_{W',W}\\
&&+\left\langle{\partial\over\partial t}b(x,\by)(u_0+\sum_{i=1}^n\nabla_{y_i}\fu_i(t,x,\by_i)\Big|_{t=s},\fw_i(s,\cdot,\cdot)\right\rangle_{V_i',V_i}\\
&&-\left\langle{\partial\over\partial t}\int_\bY b(x,\by)(u_0(t,x)+\sum_{i=1}^n\nabla_{y_i}\fu_i(t,x,\by_i))d\by\Big|_{t=0},w_0(0,\cdot)\right\rangle_{W',W}\\
&&-\left\langle{\partial\over\partial t}b(x,\by)(u_0+\sum_{i=1}^n\nabla_{y_i}\fu_i(t,x,\by_i)\Big|_{t=0},\fw_i(0,\cdot,\cdot)\right\rangle_{V_i',V_i}\\
&&=\int_0^s\int_D\int_\bY b(x,\by)\left({\partial u_0\over\partial t}+\sum_{i=1}^n{\partial\over\partial t}\nabla_{y_i}\fu_i\right)\cdot\left({\partial w_0\over\partial t}+\sum_{i=1}^n{\partial\over\partial t}\nabla_{y_i}\fw_i\right)d\by dx dt\\
&&-\int_0^s\int_D\int_\bY a(x,\by)\left(\curl u_0+\sum_{i=1}^n\scurl_{y_i} u_i\right)\cdot\left(\curl w_0+\sum_{i=1}^n\scurl_{y_i}w_i\right)d\by dxdt\\
&&=\int_0^s\frac12{d\over dt}\int_D\int_\bY b(x,\by)(u_0+\sum_{i=1}^n\nabla_{y_i}\fu_i)\cdot(u_0+\sum_{i=1}^n\nabla_{y_i}\fu_i)d\by dxdt-\\
&&\int_0^s\frac12{d\over dt}\int_D\int_\bY a(x,\by)(\curl w_0+\sum_{i=1}^n\scurl_{y_i}w_i)\cdot(\curl w_0+\sum_{i=1}^n\scurl_{y_i}w_i)d\by dxdt 
\eeqas
due to $u_0={\partial w_0\over\partial t}$, $u_i={\partial w_i\over\partial t}$ and $\fu_i={\partial\fw_i\over\partial t}$ for $i=1,\ldots,n$ from \eqref{eq:w}. Using \eqref{eq:fuiv0i} and \eqref{eq:fuivii}, we have
\beqas
&&0=\frac12\int_D\int_\bY b(x,\by)(u_0(s,x)+\sum_{i=1}^n\nabla_{y_i}\fu_i(s,x,\by_i))\cdot(u_0(s,x)+\sum_{i=1}^n\nabla_{y_i}\fu_i(s,x,\by_i))d\by dx+\\
&&\frac12\int_D\int_\bY a(x,\by)(\curl w_0(0,x)+\sum_{i=1}^n\scurl_{y_i} w_i(0,x,\by_i))\cdot(\curl w_0(0,x)+\sum_{i=1}^n\scurl_{y_i} w_i(0,x,\by_i))d\by dx.
\eeqas
We thus deduce that $u_0(s)=0$, $\nabla_{y_i}\fu_i(s)=0$ $\forall\,s$, $\curl w_0(0)=0$ and $\scurl_{y_i}w_i(0)=0$. This means that 
\beqas
\int_0^s\curl u_0(\sigma)d\sigma=0\ \ \mbox{and}\ \ \int_0^s\scurl_{y_i}u_i(\sigma)d\sigma=0
\eeqas
for all $s$. Thus for all $\sigma$, $\curl u_0(\sigma)=0$ and $\scurl_{y_i}(\sigma)=0$. 
\eproof

\section{Homogenized equation}
In this section, we use the multiscale homogenized problem \eqref{eq:msprob} to deduce the homogenized equation. From \eqref{eq:msprob}, we have that for all $\fv_n\in V_n$ and all $q\in {\cal D}(0,T)$
\[
\int_D\int_\bY b(x,\by)\left(\int_0^T(u_0+\sum_{i=1}^{n-1}\nabla_{y_i}\fu_i)q''(t)dt+\nabla_{y_n}\int_0^T\fu_nq''(t)dt\right)\cdot\nabla_{y_n}\fv_n d\by dx=0.
\]
Thus 
\[
\int_0^T \fu_n q''(t)dt=\left(\int_0^T(u_0+\sum_{i=1}^{n-1}\nabla_{y_i}\fu_i)_kq''(t)dt\right)w_n^k
\]
where $w_n^k\in L^2(D\times\bY_{n-1},H^1_\#(Y_n)/\IR)$ is the solution of the cell problem
\be
\nabla_{y_n}\cdot(b(x,\by)(e^k+\nabla_{y_n}w_n^k)=0;
\label{eq:cellbn}
\ee
$e^k$ is the $k$-th unit vector in $\IR^d$. Therefore
\[
\int_0^T\left(\nabla_{y_n}\fu_n-(u_0+\sum_{i=1}^{n-1}\nabla_{y_i}\fu_i)_k\nabla_{y_n}w^k_n\right)q''(t)dt=0
\]
so
\[
\int_0^T\left({\partial\over\partial t}\nabla_{y_n}\fu_n-{\partial\over\partial t}(u_0+\sum_{i=1}^{n-1}\nabla_{y_i}\fu_i)_k\nabla_{y_n}w^k_n\right)q'(t)dt=0.
\]
From this we have
\be
{\partial\over\partial t}\nabla_{y_n}\fu_n={\partial\over\partial t}(u_0+\sum_{i=1}^{n-1}\nabla_{y_i}\fu_i)_k\nabla_{y_n}w_n^k+G_n(x,\by_n)
\label{eq:Gn}
\ee
for a function $G_n(x,\by_n)$ in $L^2(D\times\bY_n)$. We then have $\forall\,\fv_{n-1}\in V_{n-1}$
\[
\int_0^T\int_D\int_{\bY_{n-1}}b^{n-1}(x,\by_{n-1})(u_0+\sum_{i=1}^{n-1}\nabla_{y_i}\fu_i)\cdot\nabla_{y_{n-1}}\fv_{n-1} q''(t)d\by_{n-1}dxdt=0,
\]
where $b^{n-1}(x,\by_{n-1})$ is the $(n-1)$th level homogenized coeffcient which is defined by
\be
b^{n-1}_{ij}(x,\by_{n-1})=\int_{Y_n}b_{kl}(x,\by)(\delta_{jl}+{\partial w_n^j\over\partial y_{nl}})(\delta_{ik}+{\partial w_n^i\over\partial y_{nk}})dy_n.
\label{eq:bnminus1}
\ee
Similarly we have
\[
{\partial\over\partial t}\nabla_{y_{n-1}}\fu_{n-1}={\partial\over\partial t}(u_0+\sum_{i=1}^{n-2}\nabla_{y_i}\fu_i)_k\nabla_{y_{n-1}}w^k_{n-1}+G_{n-1}(x,\by_{n-1}),
\]
where $G_{n-1}(x,\by_{n-1})\in L^2(D\times\bY_{n-1})$, and $w^k_{n-1}\in L^2(D\times\bY_{n-2}, H^1_\#(Y_{n-1})/\IR)$ satisfies the cell problem
\[
\nabla_{y_{n-1}}\cdot(b^{n-1}(e^k+\nabla_{y_{n-1}} w^k_{n-1}))=0.
\]
Recursively, letting $b^n(x,\by_n)=b(x,\by)$, we have for all $i=0,\ldots,n-1$:
\[
\int_0^T \fu_i q''(t)dt=\left(\int_0^T(u_0+\sum_{j=1}^{i-1}\nabla_{y_j}\fu_j)_kq''(t)dt\right)w^k_i,
\]
where $w^k_i\in L^2(D\times \bY_{i-1},H^1_\#(Y_i)/\IR)$ is the solution of the cell problem
\be
\nabla_{y_i}\cdot(b^i(x,\by_k)(e^k_i+\nabla_{y_i} w^k_i))=0.
\label{eq:cellbi}
\ee
From an argument as above, we have 
\[
{\partial\over\partial t}\nabla_{y_i}\fu_i={\partial\over\partial t}(u_0+\sum_{j=1}^{i-1}\nabla_{y_j}\fu_j)_k\nabla_{y_i} w^k_i+G_i(x,\by_i)
\]
for a function $G_i(x,\by_i)\in L^2(D\times\bY_i)$. 
The positive definite matrix function $b^0(x)$, which is defined as
\be
b^0_{pq}(x)=\int_Yb_{kl}^1(x,y_1)\left(\delta_{ql}+{\partial w^q_1\over\partial y_{1l}}\right)\left(\delta_{pk}+{\partial w^p_1\over\partial y_{1k}}\right)dy_1
\label{eq:b0}
\ee
is the homogenized coefficient. It satisfies
\beqas
\int_0^T\int_D\int_\bY b(x,\by)(u_0+\sum_{i=1}^n\nabla_{y_i}\fu_i)\cdot(v_0+\sum_{i=1}^n\nabla_{y_i}\fv_i)q''(t)d\by dx dt\\
=\int_0^T\int_D b^0(x) u_0\cdot v_0 q''(t)dxdt
\eeqas
for all $v_0\in W$ and $\fv_i\in V_i$. 

From \eqref{eq:msprob} 
\[
\int_D\int_\bY a(x,\by)\left(\curl u_0+\sum_{i=1}^n\scurl_{y_i}u_i\right)\cdot\scurl_{y_n}v_n d\by dx=0
\]
for all $v_n \in W_n$. For each $l=1,...,d$, let $N^l_n \in V_n$ be the solution of
\[
\scurl_{y_n}(a(x,\by)(e^l+\scurl_{y_n}N^l_n))=0. 
\]
We can write $u_n$ as
\[
u_n=N_n^l \big((\curl u_0)_l+(\scurl_{y_1}u_1)_l+...+(\scurl_{y_{n-1}}u_{n-1})_l\big).
\]
For all $v_{n-1}\in W_{n-1}$,
\[
\int_D\int_\bY a(x,\by)\left(\curl u_0+\sum_{i=1}^{n-1}\scurl_{y_i}u_i+\curl_{y_n}N_n^l \big((\curl u_0)_l+\sum_{i=1}^{n-1}(\curl_{y_i}u_i)_l\big)\right)\cdot\curl_{y_{n-1}}v_{n-1} d\by dx=0
\]
i.e. 
\[
\int_D\int_\bY a_{ij}(x,\by)\left(\delta_{lj}+(\curl_{y_n}N_n^l)_j\right)\left((\curl u_0)_l+\sum_{i=1}^{n-1}(\scurl_{y_i}u_i)_l\right)(\curl_{y_{n-1}}v_{n-1})_i d\by dx=0
\]
for all $v_{n-1} \in W_{n-1}$.
Let
\[
a^{n-1}_{pq}(x,y_1,..,y_{n-1})=\int_{Y_n}a_{pk}(x,\by)\left(\delta_{kq}+(\scurl_{y_n}N_n^q)_k\right)dy_n.
\]
We have that
\[
u_{n-1}=N_{(n-1)}^l \big((\curl u_0)_l+(\scurl_{y_1}u_1)_l+...+(\scurl_{y_{n-2}}u_{n-2})_l\big)
\]
where $N_{(n-1)}^l $ satisfies the cell problem
\[
\scurl_{y_{n-1}}(a^{n-1}(x,\by_{n-1})(e^l+\scurl_{y_{n-1}}N^l_{n-1}))=0.
\]
Letting $a^n=a$, we then have, recursively,
\[
u_{i}=N_{i}^l \big((\curl u_0)_l+(\scurl_{y_1}u_1)_l+...+(\scurl_{y_{i-1}}u_{i-1})_l\big)
\]
where $N_i^l \in  W_i$ satisfies the cell problem
\be
\scurl_{y_i}(a^i(e^l+\scurl_{y_i}N^l_i))=0,
\label{eq:cellai}
\ee
i.e. 
\[
    \int_D\int_{\bY_{i}} a^i(e^l+\scurl_{y_{i}}N^l_{i})\cdot\scurl_{y_{i}}v_{i} d\by_i=0
\]
for all $v_{i} \in W_i$. For $i=1,\ldots,n-1$, the $i$th level homogenized coefficient $a^i$ is defined as
\[
a^{i}_{pq}(x,y_1,..,y_{i})=\int_{Y_{i+1}}a^{i+1}_{pk}\left(\delta_{kq}+(\scurl_{y_{i+1}}N_{i+1}^q)_k\right)dy_{i+1}.
\]
Continuing this process, we finally get the homogenized coefficient $a^0(x)$ as
\be
a^0_{pq}(x)=\int_{Y_{1}}a^1_{pk}\left(\delta_{kq}+(\scurl_{y_{1}}N_{1}^q)_k\right).
\label{eq:a0}
\ee
The homogenization equation is
\[
\int_0^T\int_D b^0(x)u_0(x)\cdot v_0(x)q''(t)dxdt+\int_0^T\int_D a^0\curl u_0\cdot\curl v_0 q(t) dxdt=\int_0^T\int_D f(t,x)\cdot v_0(x) q(t)dxdt
\]
i.e. 
\beq
b^0(x) \frac{\partial^2u_0}{\partial t^2}(t,x)+\curl(a^0(x)\curl u_0(x))=f(t,x).
\label{eq:homogenizedeq}
\eeq
Now we derive the initial conditions. From \eqref{eq:u0i}, $u_0(0)=g_0$. As a distribution in $V_n'$, 
\[
{\partial^2\over\partial t^2}b(x,\by)(u_0+\sum_{i=1}^n\nabla_{y_i}\fu_i)=0,
\]
so for all $t$
\[
{\partial\over\partial t}b(x,\by)(u_0+\sum_{i=1}^n\nabla_{y_i}\fu_i)={\partial\over\partial t}b(x,\by)(u_0+\sum_{i=1}^n\nabla_{y_i}\fu_i)\Big|_{t=0},
\]
i.e.
\[
\int_D\int_\bY b(x,\by)\left({\partial u_0\over\partial t}+\sum_{i=1}^n{\partial\over\partial t}\nabla_{y_i}\fu_i\right)\cdot\nabla_{y_i}\fv_n d\by dx=\int_D\int_\bY b(x,\by)g_1\cdot\nabla_{y_n}\fv_n d\by dx.
\]
From \eqref{eq:Gn} we have
\beqas
&&\int_D\int_\bY b(x,\by)g_1\cdot\nabla_{y_n}\fv_n d\by dx\\
&&=\int_D\int_\bY b(x,\by)\left({\partial\over\partial t}(u_0+\sum_{i=1}^{n-1}\nabla_{y_i}\fu_i)_k(e^k+\nabla_{y_n}w^k_n)+G_n(x,\by_{n})\right)\cdot\nabla_{y_n}\fv_n (y_n)d\by dx\\
&&=\int_D\int_\bY b(x,\by)G_n(x,\by_{n})\cdot\nabla_{y_n}\fv_n(y_n)d\by dx
\eeqas
due to \eqref{eq:cellbn}. From \eqref{eq:Gn} we have
\[
\scurl_{y_n}G_n(x,\by_n)=0,\ \ \mbox{and}\ \ \int_{Y_n}G_n(x,\by_n)dy_n=0.
\]
Thus there is a function $\tilde G_n\in L^2(D\times \bY_{n-1},H^1_\#(Y_n)/\IR)$ such that $G_n(x,\by_n)=\nabla_{y_n}\tilde G_n(x,\by_n)$. From 
\[
\int_D\int_{\bY} b(x,\by)(-g_1+\nabla_{y_n}\tilde G_n(x,\by))\cdot\nabla_{y_n}\fv(y_n)d\by dx=0,
\]
we deduce that 
\be
\tilde G_n(x,\by)=-g_{1k}w^k_n\ \ \mbox{so}\ \ G_n(x,\by)=-g_{1k}\nabla_{y_n}w^k_n
\label{eq:tildeGn}
\ee
 where $w^k_n$ is the solution of the cell problem \eqref{eq:cellbn}. As a function in $V_{n-1}'$, ${\partial^2\over\partial t^2}b(x,\by)(u_0+\sum_{i=1}^n\nabla_{y_i}\fu_i)=0$ so 
\[
{\partial\over\partial t}b(x,\by)(u_0+\sum_{i=1}^n\nabla_{y_i}\fu_i)={\partial\over\partial t}b(x,\by)(u_0+\sum_{i=1}^n\nabla_{y_i}\fu_i)\Big|_{t=0}.
\]
From \eqref{eq:fuivii}, we have
\[
\int_D\int_\bY b(x,\by)\left({\partial u_0\over\partial t}+\sum_{i=1}^n{\partial\over\partial t}\nabla_{y_i}\fu_i\right)\cdot\nabla_{y_{n-1}}\fv_{n-1}d\by dx=\int_D\int_\bY b(x,\by)g_1\cdot\nabla_{y_{n-1}}\fv_{n-1}d\by dx.
\]
From \eqref{eq:Gn} and \eqref{eq:tildeGn} we have
\beqas
&&\int_D\int_\bY b(x,\by)(e^k+\nabla_{y_n}w^k_n)\left({\partial u_0\over\partial t}+\sum_{j=1}^{n-1}{\partial\over\partial t}\nabla_{y_j}\fu_j\right)_k\cdot\nabla_{y_{n-1}}\fv_{n-1}d\by dx\\
&&-\int_D\int_\bY b(x,\by)g_{1k}\nabla_{y_n}w^k_n\cdot\nabla_{y_{n-1}}\fv_{n-1}d\by dx=\int_D\int_\bY b(x,\by)g_1(x)\cdot\nabla_{y_{n-1}}\fv_{n-1}d\by dx.
\eeqas
Thus
\beqas
&&\int_D\int_{\bY_{n-1}}b^{n-1}(x,\by_{n-1})\left({\partial u_0\over\partial t}+\sum_{j=1}^{n-1}{\partial\over\partial t}\nabla_{y_j}\fu_j\right)\cdot\nabla_{y_{n-1}}\fv_{n-1}d\by_{n-1}dx\\
&&=\int_D\int_\bY b(x,\by)(e^k+\nabla_{y_n}w^k_n)g_{1k}\nabla_{y_{n-1}}\fv_{n-1}d\by dx\\
&&=\int_D\int_{\bY_{n-1}}b^{n-1}(x,\by_{n-1})g_1(x)\cdot\nabla_{y_{n-1}}\fv_{n-1}d\by_{n-1}dx.
\eeqas
Continuing this, we find that
\[
G_i(x,\by_i)=-g_{1k}\nabla_{y_i}w^k_i.
\]
Therefore, for all $t$
\beqas
&&\int_D\int_\bY b(x,\by)({\partial u_0\over\partial t}(t,x)+\sum_{i=1}^n{\partial\over\partial t}\nabla_{y_i}\fu_i(t,x,\by))\cdot v_0dxd\by\\
&&=\int_D b_0(x){\partial u_0\over\partial t}(t,x)\cdot v_0(x) dx-\sum_{i=1}^n\int_D\int_\bY b_i(x,\by_i)\nabla w^k_ig_{1k}\cdot v_0 d\by dx.
\eeqas
Therefore for $q\in C^\infty([0,T])$ with $q(T)=0$, using the fact that ${\partial\over\partial t}b(x,\by)(u_0+\sum_{i=1}^n\nabla_{y_i}\fu_i)$ is continuous as a map from $[0,T]$ to $W'$, we have
\beqas
&&\int_0^T\left\langle{\partial^2\over\partial t^2}\int_\bY b(x,\by)(u_0(t,x)+\sum_{i=1}^n\nabla_{y_i}\fu_i(t,x,\by))d\by,v_0\right\rangle_{W',W}q(t)dt\\
&&=\int_0^T{\partial\over\partial t}\left(\left\langle{\partial\over\partial t}\int_\bY b(x,\by)(u_0(t,x)+\sum_{i=1}^n\nabla_{y_i}\fu_i(t,x,\by))d\by,v_0\right\rangle_{W',W}q(t)\right) dt\\
&&\qquad\qquad-\int_0^T\left\langle\int_\bY b(x,\by)({\partial u_0\over\partial t}+\sum_{i=1}^n{\partial\over\partial t}\nabla_{y_i}\fu_i),v_0\right\rangle q'(t)dt\\
&&=-\left\langle b^0{\partial u_0\over\partial t}\Big|_{t=0}, v_0\right\rangle_{W',W} q(0)+\sum_{i=1}^n\int_D\int_\bY b_i(x,\by_i)\nabla w^k_i g_{1k}\cdot v_0 d\by dx q(0)\\
&&-\int_0^T\left\langle\int_\bY b(x,\by)({\partial u_0\over\partial t}+\sum_{i=1}^n{\partial\over\partial t}\nabla_{y_i}\fu_i),v_0\right\rangle q'(t)dt.
\eeqas
From \eqref{eq:sec2i1} and \eqref{eq:expr},  we deduce
\beqas
-\left\langle b^0{\partial u_0\over\partial t}\Big|_{t=0}, v_0\right\rangle q(0)+\sum_{i=1}^n\int_D\int_\bY b_i(x,\by_i)\nabla w^k_i g_{1k}\cdot v_0 d\by dx q(0)
=-\int_D\int_\bY b(x,\by)g_1(x)\cdot v_0(x)dx q(0).
\eeqas
We then have
\beqas
\left\langle b^0{\partial u_0\over\partial t}\Big|_{t=0}, v_0\right\rangle
=\int_D\int_\bY b(x,\by)g_1(x)\cdot v_0(x)dx+\sum_{i=1}^n\int_D\int_\bY b_i(x,\by_i)\nabla w^k_i g_{1k}\cdot v_0 d\by dx.
\eeqas
We note that
\beqas
\int_D\int_\bY b(x,\by)g_1(x)\cdot v_0(x)dx +\int_D\int_\bY b(x,\by)\nabla_{y_n} w^k_ng_{1k}\cdot v_0d\by dx
=\int_D\int_\bY b(x,\by)(e^k+\nabla_{y_n}w^k_n)g_{1k}(x)\cdot v_0d\by dx\\
=\int_D\int_{\bY_{n-1}} b^{n-1}(x,\by_{n-1})g_1\cdot v_0 d\by_{n-1} dx.
\eeqas
Continuing this we have
\[
\int_D\int_\bY b(x,\by)g_1(x)\cdot v_0(x)dx+\sum_{i=1}^n\int_D\int_\bY b_i(x,\by_i)\nabla w^k_i g_{1k} d\by dx=\int_D b^0(x)g_1\cdot v_0dx
\]
Thus as distribution in $V'$
\be
b^0{\partial u_0\over\partial t}\Big|_{t=0}=b^0g_1.
\label{eq:du0dt}
\ee
Therefore, $u^0$ is the solution of the problem \eqref{eq:homogenizedeq} with the initial condition $u^0(0)=g^0$ and \eqref{eq:du0dt} which has a unique solution.

The solution $\bu$ is written in terms of $u_0$ as
\be
u_i=N_i^{r_{i-1}}(\delta_{r_{i-1}r_{i-2}}+(\scurl_{y_{i-1}}N^{r_{i-2}}_{r_{i-1}}))\ldots(\curl u_0)_{r_0},
\label{eq:ui}
\ee
and
\be
{\partial\over\partial t}\nabla_{y_i}\fu_i={\partial u_{0r_0}\over\partial t}(\delta_{r_0r_1}+{\partial w_1^{r_0}\over\partial y_{1r_1}})(\delta_{r_1r_2}+{\partial w_2^{r_1}\over\partial  y_{2r_2}})\ldots\nabla_{y_i}w^{r_{i-1}}_i-g_{1r_0}(\delta_{r_0r_1}+{\partial w_1^{r_0}\over\partial y_{1r_1}})(\delta_{r_1r_2}+{\partial w_2^{r_1}\over\partial  y_{2r_2}})\ldots\nabla_{y_i}w^{r_{i-1}}_i.
\label{eq:dtfui}
\ee
Given that at $t=0$, $\nabla_{y_i}\fu_i=0$, we then have
\begin{eqnarray}
\nabla_{y_i}\fu_i=u_{0r_0}(\delta_{r_0r_1}+{\partial w_1^{r_0}\over\partial y_{1r_1}})(\delta_{r_1r_2}+{\partial w_2^{r_1}\over\partial  y_{2r_2}})\ldots\nabla_{y_i}w^{r_{i-1}}_i-g_{1r_0}(\delta_{r_0r_1}+{\partial w_1^{r_0}\over\partial y_{1r_1}})(\delta_{r_1r_2}+{\partial w_2^{r_1}\over\partial  y_{2r_2}})\ldots\nabla_{y_i}w^{r_{i-1}}_it\nonumber\\
-g_{0r_0}(\delta_{r_0r_1}+{\partial w_1^{r_0}\over\partial y_{1r_1}})(\delta_{r_1r_2}+{\partial w_2^{r_1}\over\partial  y_{2r_2}})\ldots\nabla_{y_i}w^{r_{i-1}}_i.
\label{eq:nablafui}
\end{eqnarray}
Without loss of generality, we let
\begin{eqnarray}
\fu_i=u_{0r_0}(\delta_{r_0r_1}+{\partial w_1^{r_0}\over\partial y_{1r_1}})(\delta_{r_1r_2}+{\partial w_2^{r_1}\over\partial  y_{2r_2}})\ldots w^{r_{i-1}}_i-g_{1r_0}(\delta_{r_0r_1}+{\partial w_1^{r_0}\over\partial y_{1r_1}})(\delta_{r_1r_2}+{\partial w_2^{r_1}\over\partial  y_{2r_2}})\ldots w^{r_{i-1}}_it\nonumber\\
-g_{0r_0}(\delta_{r_0r_1}+{\partial w_1^{r_0}\over\partial y_{1r_1}})(\delta_{r_1r_2}+{\partial w_2^{r_1}\over\partial  y_{2r_2}})\ldots w^{r_{i-1}}_i.
\label{eq:fui}
\end{eqnarray}

\section{Regularity of the solution}
%
To deduce the homogenization errors in the next section, we need regularity for the solution $u^0$ of the homogenized equation \eqref{eq:homogenizedeq}, and of the solutions of the cell problems \eqref{eq:cellbi} and \eqref{eq:cellai}. We assume:
\begin{assumption}
\label{assum:regularityab}
The matrix functions $a$ and $b$ belong to $C^1(\bar D,C^2(\bar Y_1,\ldots,C^2(\bar Y_n)\ldots))^{d\times d}$.
\end{assumption}
With this assumption, we  have
\begin{proposition}
\label{prop:regularityNomega}
Under Assumption \ref{assum:regularityab}, for all $i,r=1,\ldots,d$, $\scurl_{y_i}N_i^r\in C^1(\bar D,C^2(\bar Y_1,\ldots,C^2(\bar Y_{i-1},H^2(Y_i))\ldots))$ and $w_i^r\in C^1(\bar D,C^2(\bar Y_1,\ldots,C^2(\bar Y_{i-1},H^3(Y_i))\ldots))$.
\end{proposition}
We refer to  \cite{Tiep1} for a proof of this proposition.

We have the following regularity results for the solution $u_0$ of the homogenized equation \eqref{eq:homogenizedeq}. 

\begin{proposition}
Under Assumption \ref{assum:regularityab}, assume 
\beq
\begin{cases}
f \in H^2(0,T; H),&\\
g_1 \in W,&\\
(b^0)^{-1}[f(0)-\curl (a^0(x)\curl g_0)] \in W,&\\
(b^0)^{-1}[\frac{\partial f}{\partial t}(0)-\curl (a^0(x)\curl g_1)] \in H,&
\end{cases}
\label{eq:Compatibility1}
\eeq
then 
\beq
\frac{\partial^2 u_0}{\partial t^2}\in L^\infty(0,T; W),\ \ \frac{\partial^3 u_0}{\partial t^3}\in L^\infty(0,T; H),\ \ \mbox{and}\ \ {\partial^3\over\partial t^3}\nabla_{y_i}\fu_i\in L^\infty(0,T;L^2(D\times\bY)).
\label{eq:reg1}
\eeq

Further, if
\beq
\begin{cases}
f \in H^3(0,T; H),&\\
g_1 \in W,&\\
(b_0)^{-1}[f(0)-\curl (a^0(x)\curl g_0)] \in W,&\\
(b_0)^{-1}[\frac{\partial f}{\partial t}(0)-\curl (a^0(x)\curl g_1)] \in W,&\\
(b_0)^{-1}[\frac{\partial^2 f}{\partial t^2}(0)-\curl (a^0(x)\curl ((b_0)^{-1}(f_0-\curl (a^0(x)\curl g_0)))) \in H,&
\end{cases}
\label{eq:Compatibility2}
\eeq
then 
\beq
\frac{\partial^3 u_0}{\partial t^3}\in L^\infty(0,T; V),\ \ \frac{\partial^4 u_0}{\partial t^4}\in L^\infty(0,T; H),\ \ \mbox{and}\ {\partial^4\over\partial t^4}\nabla_{y_i}\fu_i\in L^\infty(0,T;L^2(0,T;L^2(D\times\bY)).  
\label{eq:reg2}
\eeq
\end{proposition}

\bproof
We use the regularity theory of general hyperbolic equations (see, e.g., Wloka \cite{Wloka}, Chapter 5). 
From \eqref{eq:Compatibility1} we have that
\be
b^0\frac{\partial ^2}{\partial t ^2}\left (\frac{\partial u_0}{\partial t }\right)+\curl \left(a^0\curl \frac{\partial u_0}{\partial t }\right)=\frac{\partial f}{\partial t }
\label{eq:e1}
\ee
with compatibility initial conditions
\[
\frac{\partial u_0}{\partial t }(0)=g_1 \in W,~~\frac{\partial }{\partial t }\frac{\partial u_0}{\partial t }(0)=(b^0)^{-1}[f(0)-\curl (a^0\curl g_0) \in W
\]
and
\be
\label{eq:e2}
b^0\frac{\partial ^2}{\partial t ^2}\left(\frac{\partial^2 u_0}{\partial t^2 }\right)+\curl \left(a^0\curl \frac{\partial^2 u_0}{\partial t^2 }\right)=\frac{\partial^2 f}{\partial t^2 },
\ee
with compatibility initial conditions
\[
\frac{\partial^2 u_0}{\partial t^2 }(0)=(b^0)^{-1}[f(0)-\curl \left(a^0\curl g_0\right) \in W
\ \ 
\mbox{and}
\ \ 
\frac{\partial }{\partial t }\frac{\partial^2 u_0}{\partial t^2 }(0)=(b^0)^{-1}[\frac{\partial f}{\partial t}(0)-\curl (a^0\curl g_1)] \in  H.
\]
We thus deduce that
\[
{\partial^2u_0\over\partial t^2}\in L^\infty(0,T;W)\ \ \mbox{and}\ \ {\partial^3u_0\over\partial t^3}\in L^\infty(0,T;H).
\]
From \eqref{eq:nablafui} and  Proposition \ref{prop:regularityNomega}, we deduce that
\[
{\partial^3\over\partial t^3}\nabla_{y_i}\fu_i\in L^\infty(0,T;L^2(D\times\bY)).
\]
Similarly, we deduce regularity \eqref{eq:reg2} from \eqref{eq:Compatibility2}. \eproof

For the regularity of $u_0$, we have the following result.
\bp\label{prop:regularityu0}
Under Assumption \ref{assum:regularityab}, if $D$ is a Lipschitz polygonal domain, $f\in H^1(0,T;H)$, $g_0\in H^1(\curl,D)$ and $g_1\in W$, ${\rm div}f\in L^\infty(0,T;L^2(D))$, ${\rm div}(b^0g_0)\in L^2(D)$ and ${\rm div}(b^0g_1)\in L^2(D)$,  there is a constant $s\in(0,1]$ such that $u_0\in L^\infty(0,T;H^s(\curl,D))$.
\epr
\bproof
Using Proposition \ref{prop:regularityNomega}, equations \eqref{eq:a0} and \eqref{eq:b0}, we have that $a^0, b^0\in C^1(\bar D)^{d\times d}$.
As $f\in H^1(0,T;H)$ and $g_0\in H^1(\curl,D)$, we have that $(b^0)^{-1}[f-\curl(a^0\curl g^0)]\in H$. The compatibility initial conditions hold so that ${\partial^2u_0\over\partial t^2}\in L^\infty(0,T;H)$. Thus
\[
\curl(a^0\curl u_0)=f-b^0{\partial^2u_0\over\partial t^2}\in L^\infty(0,T;H).
\]
Let $U(t)=a^0\curl u_0(t)$. As ${\rm div}((a^0)^{-1}U(t))=0$ and $(a^0)^{-1}U(t)\cdot\nu=0$, there is a constant $c$ and a constant $s\in (0,1]$ which depend on $a^0$ and the domain $D$ so that
\[
\|U(t)\|_{H^s(D)^3}\le c(\|\curl U(t)\|_{L^2(D)^3}+\|U(t)\|_{L^2(D)^3})
\]
so $U\in L^\infty(0,T;H^s(D)^3)$. As $\curl u_0(t)=(a^0)^{-1}U(t)$ and $(a^0)^{-1}\in C^1(\bar D)^{d\times d}$, $\curl u_0\in L^\infty(0,T;H^s(D))$. 
We note that
\[
{\rm div}\left(b^0{\partial^2u_0\over\partial t^2}\right)={\rm div} f,
\]
so 
\[
{\rm div}(b^0u_0(t))=\int_0^t\int_0^s{\rm div}f(s)ds+t{\rm div}(b^0g_1)+{\rm div}(b^0 g_0)\in L^\infty(0,T;L^2(D)).
\]
From Theorem 4.1 of Hiptmair \cite{Hiptmair2002}, we deduce that there is a constant $s\in (0,1]$ (we take it as the same constant as above), so that 
\[
\|u_0(t)\|_{H^s(D)^3}\le c(\|u(t)\|_{H(\curl,D)}+\|{\rm div}(b^0u_0(t))\|_{L^2(D)^3}).
\]
Thus $u_0\in L^\infty(0,T;H^s(\curl,D))$.   
\eproof

Similarly, we can deduce the regularity for ${\partial^2u_0\over\partial t^2}$.
\bp\label{prop:regularityd2u0dt2}
Under Assumption \ref{assum:regularityab}, if $D$ is a Lipschitz polygonal domain, if the comparibility conditions \eqref{eq:Compatibility2} hold, and if ${\rm div}f\in L^\infty(0,T;L^2(D))$, then  there is a constant $s\in (0,1]$ such that ${\partial^2u_0\over\partial t^2}\in L^\infty(0,T;H^s(\curl,D))$.
\epr
\bproof
From equation \eqref{eq:e2}, we have
\[
\curl\left(a^0\curl{\partial^2u_0\over\partial t^2}\right)={\partial^2 f\over\partial t^2}-b^0{\partial^4u_0\over\partial t^4}\in L^\infty(0,T;H)
\]
as ${\partial^4u_0\over\partial t^4}\in L^\infty(0,T;H)$ due to \eqref{eq:Compatibility2}. Following a similar argument as in the proof of Proposition \ref{prop:regularityu0} we deduce that $\curl{\partial^2u_0\over\partial t^2}\in L^\infty(0,T;H^s(D)^3)$. We note that 
\[
{\rm div}b^0{\partial^2 u_0\over\partial t^2}={\rm div} f\in L^\infty(0,T;L^2(D)).
\]
From Theorem 4.1 of \cite{Hiptmair2002}, we deduce that ${\partial^2u_0\over\partial t^2}\in L^\infty(0,T;H^s(D)^3)$.  
\eproof

\section{Corrector for the homogenization problem}
We derive homogenization correctors in the section. For two scale problems, with sufficient regularity for the solution of the homogenized equation and the cell problems, we derive explicit homogenization errors in terms of the microscopic scales. We consider both the cases where $u_0$ and ${\partial\over\partial t}u_0$ 
belong to $L^\infty((0,T);H^1(\curl,D))$ and where they only belong to a weaker regularity space $L^\infty((0,T);H^s(\curl,D))$ for $0<s<1$, which is normally the case in polygonal domains $D$. For multiscale problems, we are not able to prove an explicit homogenization error but a general corrector is derived. 
\subsection{Corrector for two-scale problem}
%
%
For conciseness, we denote the coefficients $a(x,\by)$ and $b(x,\by)$ as $a(x,y)$ and $b(x,y)$. The cell problems become
\be
\scurl_y(a(x,y)(e^r+\scurl_yN^r(x,y))=0,
\label{eq:2sNcell}
\ee
and 
\be
\nabla_yb(x,y)\cdot(e^r+\nabla_y \omega^r(x,y))=0.
\label{eq:2swcell}
\ee
The homogenized coefficient is determined by
\be
a^0_{pq}(x)=\int_Ya_{pk}(x,y)(\delta_{kq}+(\scurl_yN^q)_k)dy
\label{eq:2sa0}
\ee
and 
\be
b^0_{pq}(x)=\int_Yb_{pl}\left(\delta_{ql}+{\partial\omega^q\over\partial y_l}\right)dy.
\label{eq:2sb0}
\ee

We then have the following
\begin{proposition}\label{prop:2sregularhomerror}
For two-scale problems, assume that
$g_0=0$, $g_1 \in H^1(D)^3\bigcap W$, 
$f \in H^1 (0,T; H)$, 
$u_0 \in L^\infty(0,T;H^1(\curl; D))$, ${\partial u_0\over\partial t}\in L^\infty(0,T;H^1(\curl,D))$, ${\partial^2u_0\over\partial t^2}\in L^\infty(0,T;H(\curl,D)\bigcap H^1(D)^3)$, $N^r \in C^1(\bar{D}, C(\bar{Y}))^3$, $\scurl_y N^r \in  C^1(\bar{D}, C(\bar{Y}))^3 $ and $\omega^r \in C^1(\bar{D}, C(\bar{Y}))$ for all $r=1,2,3$. There exists a constant $c$ that does not depend on $\varepsilon$ such that
\[
\left\|\frac{\partial u^\varepsilon}{\partial t}-\left[\frac{\partial u_0}{\partial t}+\nabla_y\frac{\partial \frak u_1}{\partial t}\left(\cdot, \cdot, \frac{\cdot}{\varepsilon} \right ) \right]\right\|_{L^{\infty}(0,T; H)}+\left\|\curl u^\varepsilon-\left[\curl u_0+ \curl_y u_1 \left(\cdot, \cdot, \frac{\cdot}{\varepsilon} \right )\right]\right\|_{L^{\infty}(0,T;H)} \leq c \varepsilon^{1/2}.
\]

\end{proposition}

\bproof 
We note that ${\partial\ue\over\partial t}$ satisfies
\be
{\partial^2\over\partial t^2}\left({\partial\ue\over\partial t}\right)+\curl\left(a^\ep\curl{\partial\ue\over\partial t}\right)={\partial f\over\partial t}
\label{eq:duepdt}
\ee
with the initial condition ${\partial\ue\over\partial t}(0)=g_1\in W$ and ${\partial^2\ue\over\partial t^2}(0)=f(0)\in H$ (due to $g_0=0$). We therefore deduce that
\[
\left\|{\partial\ue\over\partial t}\right\|_{L^\infty(0,T;W)}\le c\left(\left\|{\partial f\over\partial t}\right\|_{L^2(0,T;H)}+\|g_1\|_V+\|f(0)\|_H\right)
\]
where $c$ only depends on the constants $\alpha$ and $\beta$ in \eqref{eq:coercive}. Thus ${\partial\ue\over\partial t}$ is uniformly bounded in $L^\infty(0,T;W)$ for all $\ep$. 

We consider the function 
\[
u_1^\ep(t,x)=u_0(t,x)+\ep N^r\left(x,{x\over\ep}\right)\left(\curl u_{0}(t,x)\right)_r+\ep\nabla\frak u_1\left(t,x,{x\over\ep}\right)
\] 
where we have from \eqref{eq:fui}
\[
\frak u_1(t,x,\frac{x}{\ep})= w^r(x, \frac{x}{\ep})\left(u_{0r}(t,x)-g_{1r}(x)t-g_{0r}(x)\right).
\]
We first show that
\[
\left\|\curl(a^\ep\curl u_1^\ep)-\curl (a^0(x)\curl u_0)\right\|_{L^\infty(0,T;W')} \leq c \ep.
\]
We have
\begin{align}
    \curl u_1^\ep
    =\curl u_0+\ep \curl _xN^r\left(x,{x\over\ep}\right)(\curl u_{0}(t,x))_r+\ep\nabla (\curl u_{0}(t,x))_r \times N^r\left(x,{x\over\ep}\right)\nonumber\\
+\curl_y N^r\left(x,{x\over\ep}\right)(\curl u_{0}(t,x))_r
\label{eq:curlu1ep}
\end{align}
Thus
\begin{align*}
a^\ep(x)\curl& u_1^\ep(t,x)- a^0\curl u_0(t,x)\\
=&a^\ep\curl u_0+a^\ep\curl_y N^r\left(x,{x\over\ep}\right)(\curl u_{0}(t,x)_r) - a^0(x)\curl u_0(t,x)\\
&+\ep a^\ep\curl _xN^r\left(x,{x\over\ep}\right)(\curl u_{0}(t,x)_r+\ep a^\ep\nabla (\curl u_{0}(t,x))_r \times N^r\left(x,{x\over\ep}\right)\\
=&G_r(x,{x\over\ep})(\curl u_0(x))_r+\ep\curl I(t,x)
\end{align*}
where the vector functions $G_r(x,y)$ are defined by
\be
(G_{r})_{i}(x,y)=a_{ir}(x,y)+a_{ij}(x,y)\left(\scurl_yN^r(x,y)\right)_j-a^0_{ir}(x),
\label{eq:Gr}
\ee
and
\beqas
I(t,x)=a\left(x,{x\over\ep}\right)\Big[\scurl_x N^r\left(x,{x\over\ep}\right)(\curl u_0(t,x))_r+\nabla(\curl u_0(t,x))_r\times N^r\left(x,{x\over\ep}\right)\Big]
\eeqas
From \eqref{eq:2sNcell}, we have that $\scurl_yG_r(x,y)=0$. Further from \eqref{eq:2sa0}, $\int_YG_r(x,y)dy=0$. From these we deduce that there are  functions $\tilde G_r(x,y)$ such that $G_r(x,y)=\nabla_y \tilde G_r(x,y)$. Since $\nabla_y\tilde G_r(x,\cdot)=G_r(x,\cdot)\in H^1(Y)^3$, we deduce that $\Delta_y\tilde G_r(x,\cdot)\in L^2(Y)$. From elliptic regularity,  $\tilde G_r(x,\cdot)\in H^2(Y)\subset C(\bar Y)$. As $G_r(x,\cdot)\in C^1(\bar D,H^1_\#(Y)^3)$,  $\tilde G_r(x,y)\in C^1(\bar D, H^2(Y))\subset C^1(\bar D,C(\bar Y))$. 
We note that 
\[
G_r(x,{x\over\ep})=\nabla_y\tilde G_r(x,{x\over\ep})=\ep\nabla\tilde G_r(x,{x\over\ep})-\ep\nabla_x\tilde G(x,{x\over\ep})
\]
Therefore, for all $\phi\in {\cal D}(D)^d$ we have
\begin{multline*}
\langle\curl(a^\ep\curl u_1^\ep)(t)-\curl(a^0\curl u_0)(t),\phi\rangle=\int_DG_r(x,{x\over\ep})(\curl u_0(t,x))_r\curl\phi(x)dx\\+\ep\int_DI(t,x)\curl\phi(x)dx=
-\ep\int_D\tilde G_r\left(x,{x\over\ep}\right){\rm div}((\curl u_0(t,x))_r\curl\phi(x))dx\\-\ep\int_D\nabla_x\tilde G_r\left(x,{x\over\ep}\right)(\curl u_0(t,x))_r\curl\phi(t,x)dx+\ep\int_DI(t,x)\curl\phi(x)dx.
\end{multline*}
We note that
\[
\int_D\tilde G_r\left(x,{x\over\ep}\right){\rm div}((\curl u_0(t,x))_r\curl\phi(t,x))dx=\int_D\tilde G_r\left(x,{x\over\ep}\right)\nabla(\curl u_0(t,x))_r\cdot\curl\phi(t,x).
\]
As  $N^j\in  C^1(\bar D,C(\bar Y))$ and  $\curl u_0\in L^\infty(0,T;H^1(D))$ we deduce that
\[
\left|\int_D I(t,x)\cdot\curl\phi(x) dx\right|\le c\|\curl\phi\|_H
\]
where $c$ is independent of $t$. From these we  conclude that 
\[
|\langle\curl(a^\ep\curl u_1^\ep)(t)-\curl(a^0(x)\curl u_0)(t),\phi\rangle|\le c\ep\|\curl\phi\|_{L^2(D)^d}.
\]
Using a density argument, we have that this holds for all $\phi\in W$. Thus
\begin{equation}
\|\curl(a^\ep\curl u_1^\ep)-\curl a^0\curl u_0\|_{L^\infty (0, T; W')}\le c\ep. \label{eq:C2}
\end{equation}
Let $\tau^\ep(x)$ be a function in ${\cal D}(D)$ such that $\tau^\ep(x)=1$ outside an $\ep$ neighbourhood of $\partial D$ and $\sup\limits_{x\in D}\ep|\nabla\tau^\ep(x)|<c$ where $c$ is independent of $\ep$. Let
\begin{equation}
w_1^\ep(t,x)=u_0(t,x)+\ep\tau^\ep(x)N^r\left(x,{x\over\ep}\right)\left(\curl u_0(t,x)\right)_r+\ep\nabla\left[\tau^\ep(x){\frak u}_1\left(t,x,{x\over\ep}\right)\right].\label{eq:C3}
\end{equation}
The function $w_1^\ep(x)$ belongs to $L^2(0,T;H_0(\curl,D))$. We note that
\begin{equation}
u_1^\ep-w_1^\ep=\ep(1-\tau^\ep(x))N^r\left(x,{x\over\ep}\right)(\curl u_0(t,x))_r+\ep\nabla\left[(1-\tau^\ep(x))w^r(x, \frac{x}{\ep})\left(u_{0r}(t,x)-g_{1r}(x)t-g_{0r}(x)\right)\right].\label{eq:C4}
\end{equation}
From this,
\begin{multline}
\curl(u_1^\ep-w_1^\ep)=\ep\curl_xN^r\left(x,{x\over\ep}\right)(\curl u_0(t,x))_r(1-\tau^\ep(x))+\curl_yN^r\left(x,{x\over\ep}\right)(\curl u_0(t,x))_r(1-\tau^\ep(x))-\\
\ep(\curl u_0(t,x))_r\nabla\tau^\ep(x)\times N^r\left(x,{x\over\ep}\right)+\ep(1-\tau^\ep(x))\nabla(\curl u_0(t,x))_r\times N^r\left(x,{x\over\ep}\right).\label{eq:C5}
\end{multline}
Since $\frac{\partial^2 u_0}{\partial t^2} \in L^\infty (0, T; H^1(D)^3)$, we have
 \begin{align*}
     \frac{\partial^2 w_1^\ep}{\partial  t^2}=&\frac{\partial^2 u_0}{\partial  t^2}+\ep \tau^\ep N^r\left(x,{x\over\ep}\right) \frac{\partial^2}{\partial t^2} \curl u_{0}(t,x)_r+ \ep \nabla \tau^\ep w^r\left(x,{x\over\ep}\right) \frac{\partial ^2  u_{0r}}{\partial t^2}(t,x)\\
     &+\ep\tau^\ep \nabla_x  w^r\left(x,{x\over\ep}\right) \frac{\partial ^2  u_{0r}}{\partial t^2}(t,x)+\tau^\ep \nabla_y w^r\left(x,{x\over\ep}\right) \frac{\partial ^2  u_{0r}}{\partial t^2}(t,x)+\ep \tau^\ep w^r\left(x,{x\over\ep}\right) \frac{\partial ^2  }{\partial t^2} \nabla u_{0r}(t,x).
 \end{align*}
 Therefore
 \begin{align}
    &b^0 \frac{\partial^2 u_0 }{\partial t^2}- b^\ep\frac{\partial^2 w_1^\ep}{\partial  t^2}\\
   &= b^0 \frac{\partial^2 u_0 }{\partial t^2}-b^\ep\frac{\partial^2 u_0}{\partial  t^2}-b^\ep \nabla_y w^r\left(x,{x\over\ep}\right) \frac{\partial ^2  u_{0r}}{\partial t^2}(t,x)\nonumber\\ 
    &-b^\ep\ep \tau^\ep N^r\left(x,{x\over\ep}\right) \frac{\partial ^2}{\partial t^2} \curl u_{0r}(t,x)- b^\ep\ep \nabla \tau^\ep w^r(x, \frac{x}{\ep}) \frac{\partial ^2  u_{0r}}{\partial t^2}(t,x)\nonumber\\
   &-b^\ep\ep\tau^\ep \nabla_x  w^r\left(x,{x\over\ep}\right) \frac{\partial ^2  u_{0r}}{\partial t^2}(t,x)+b^\ep(1-\tau^\ep) \nabla_y w^r\left(x,{x\over\ep}\right) \frac{\partial ^2  u_{0r}}{\partial t^2}(t,x)-b^\ep\ep \tau^\ep w^r\left(x,{x\over\ep}\right) \frac{\partial ^2  }{\partial t^2} \nabla u_{0r}(t,x). 
\label{eq:D8}
 \end{align}
For $\phi \in W$, we have
 \begin{align*}
    &\left\langle b^0 \frac{\partial^2 u_0 }{\partial t^2}-b^\ep\frac{\partial^2 u_0}{\partial  t^2}-b^\ep\nabla_y w^r\left(x,{x\over\ep}\right) \frac{\partial ^2  u_{0r}}{\partial t^2}(t,x), \phi\right\rangle_{H}\\
     &=\int_D\left( b^0_{ij}\frac{\partial^2 u_{0 j}}{\partial t^2}-b^\ep_{ij}\frac{\partial^2 u_{0 j}}{\partial t^2}-b^\ep_{ij}\frac{\partial w^r}{\partial y_j}\left(x,{x\over\ep}\right)\frac{\partial^2 u_{0 r}}{\partial t^2}(t,x)\right)\phi_idx\\
     &=\int_D\left( b^0_{ij}\frac{\partial^2 u_{0 j}}{\partial t^2}-b^\ep_{ij}\left(\delta_{jr}+\frac{\partial w^r}{\partial y_j}\left(x,{x\over\ep}\right)\right)\frac{\partial^2 u_{0 r}}{\partial t^2}(t,x)\right)\phi_idx\\
&=\int_D\left(b^0_{ir}(x)-b_{ir}\left(x,{x\over\ep}\right)-b_{ij}\left(x,{x\over\ep}\right){\partial w^r\over\partial y_j}(x,{x\over\ep})\right){\partial^2 u_{0r}\over\partial t^2}(t,x)\phi_i dx\\
&=\int_D F_r(x,{x\over\ep}){\partial^2u_{0r}\over\partial t^2}(t,x)\cdot\phi(x) dx
 \end{align*}
where the vector function $F_r$ is defined by 
\be
(F_r)_i(x,y)=-b_{ir}(x,y)-b_{ij}(x,y)\frac{\partial w^r}{\partial y_j}(x,y)+b^0_{ir}(x).
\label{eq:F}
\ee
Since ${\rm div}_y F_r(x,y)=0$, and $\int_Y F_r(x,y)dy=0$, there is a function $\tilde F_r(x,y)$ such that 
\be
F_r(x,y)=\curl_y\tilde F_r(x,y).
\label{eq:F1}
\ee
Detailed construction of the function $\tilde F_r$ is presented in \cite{JKO} which implies that $F_r\in C^1(\bar D,C(\bar Y))^3$. 
We note that
\[
F_r(x,{x\over\ep})=\scurl_yF_r(x,{x\over\ep})=\ep\curl\tilde F_r(x,{x\over\ep})-\ep\curl_x \tilde F_r(x,{x\over\ep}).
\]
Therefore 
 \begin{align*}
&\int_D F_r(x,{x\over\ep}){\partial^2u_{0r}\over\partial t^2}\cdot\phi(x)dx\\
&\qquad=\left|\ep \int_D \tilde F_r(x,\frac{x}{\ep})\cdot \curl \big( \frac{\partial ^2 u_{0r}}{\partial t^2}  \phi \big) dx-\ep \int_D \curl_x \tilde F_r(x,\frac{x}{\ep})\frac{\partial ^2 u_{0r}}{\partial t^2} \cdot \phi dx\right|\\
&\qquad=\left|\ep \int_D \tilde F_r(x,\frac{x}{\ep})\cdot\big(\curl \phi \frac{\partial ^2 u_{0r}}{\partial t^2}+ \phi \times \nabla \frac{\partial ^2 u_{0r}}{\partial t^2}\big)dx-\ep \int_D \curl_x \tilde F_r(x,\frac{x}{\ep})\frac{\partial ^2 u_{0r}}{\partial t^2} \cdot \phi dx\right|\\
&\qquad\le c\ep(\|\curl \phi\|_{(L^2(D))^3}+ \|\phi\|_{(L^2(D))^3})
 \end{align*}
due to the conditions ${\partial^2u_0\over\partial t^2}\in L^\infty(0,T;H^1(D)^d)$. 
Let $D^\ep \subset D$ be the $\ep$ neighbourhood of the boundary $\partial D$. We know that
\begin{equation}
\|\phi\|^2 _{L^2(D^\ep)} \leq  c \ep^2 \|\phi\|^2_{H^1(D)}+c \ep \|\phi\|^2_{L^2 (\partial D)} \leq  c \ep \|\phi\|^2_{H^1(D)} \label{eq:C7}
\end{equation}
for $\phi\in H^1(D)$ (see Hoang and Schwab \cite{HSmultirandom}). From the condition $\curl u_0\in L^\infty(0,T;H^1(D))$, ${\partial u_0\over\partial t}\in L^\infty(0,T;H^1(\curl;D))$ and ${\partial^2u_0\over\partial t^2}\in L^\infty(0,T;H^1(D))$ we deduce that: 
\be
\left\|{\partial\over\partial t}\curl u_0(t)\right\|_{L^2(D^\ep)^3}\le c\ep^{1/2},\ \ \|\curl u_0(x)_r\|_{(L^2(D^\ep))^3} \leq c \ep^{1/2}\ \ \mbox{and}\ \ \left\|\frac{\partial ^2}{\partial t^2}u_{0r}\right\|_{(L^2(D^\ep))^3} \leq c \ep^{1/2}.
\label{eq:norminDep}
\ee
From this and \eqref{eq:D8}, we get
\be
\left|\left\langle b^0 \frac{\partial^2 u_0 }{\partial t^2}- b^\ep\frac{\partial^2 w_1^\ep}{\partial  t^2} , \phi\right\rangle_H\right|\le c\ep\|\curl \phi\|_{(L^2(D))^3}+ c\ep^{1/2} \|\phi\|_{(L^2(D))^3}.
\label{eq:X1}
\ee
From
\[
  b^\ep \frac{\partial^2 u^\ep }{\partial t^2}+\curl (a^\ep \curl u^\ep)=b^0 \frac{\partial^2 u_0 }{\partial t^2}+\curl (a^0 \curl u_0)
\]
we have
 \begin{align*}
     b^\ep \frac{\partial (u^\ep -w_1^\ep)}{\partial t^2}&+\curl (a^\ep \curl (u^\ep -w_1^\ep))\\
     &=b^0 \frac{\partial^2 u_0 }{\partial t^2}- b^\ep \frac{\partial w_1^\ep }{\partial t^2}+\curl (a^\ep\curl(u_1^\ep-w_1^\ep))+\curl (a^0\curl u_0)-\curl (a^\ep\curl u_1^\ep)\\
\end{align*}
We note that
\beqas
&&{\partial w^\ep_1\over\partial t}={\partial u_0\over\partial t}+\ep\tau^\ep(x)N^r\left(x,{x\over\ep}\right){\partial\over\partial t}(\curl u_0)_r\\
&&\qquad\ \ +\left(\ep\nabla\tau^\ep(x)w^r\left(x,{x\over\ep}\right)+\ep\tau^\ep(x)\nabla_xw^r\left(x,{x\over\ep}\right)+\tau^\ep\nabla_yw^r\left(x,{x\over\ep}\right)\right)\left({\partial\over\partial t}u_{0r}(t,x)-g_{1r}(x)\right)\\
&&\qquad\ \ +\ep\tau^\ep w^r\left(x,{x\over\ep}\right)\left({\partial\over\partial t}\nabla u_{0r}(t,x)-\nabla g_{1r}(x)\right)
\eeqas
and 
\beqas
\curl{\partial w^\ep_1\over\partial t}=\curl{\partial u_0\over\partial t}+\tau^\ep(x)\left(\ep\scurl_x N^r\left(x,{x\over\ep}\right)+\scurl_yN^r\left(x,{x\over\ep}\right)\right){\partial\over\partial t}(\curl u_0)_r\\
+\left(\ep\nabla\tau^\ep(x){\partial\over\partial t}(\curl u_0)_r+\ep\tau^\ep(x)\nabla{\partial\over\partial t}(\curl u_0)_r\right)\times N^r\left(x,{x\over\ep}\right).
\eeqas
As ${\partial u_0\over\partial t}\in L^\infty(0,T;W\bigcap H^1(\curl,D))$, we deduce that ${\partial w^\ep_1\over\partial t}$ is uniformly bounded in $L^\infty(0,T;W)$. 
Since ${\partial\ue\over\partial t}(t)$ and ${\partial w_1^\ep\over\partial t}(t)$ belong to $W$ and are uniformly bounded in the norm of $W$, from \eqref{eq:X1}, we deduce that
\beqas
&&\int_0^t\left\langle b^0{\partial^2u_0\over\partial t^2}-b^\ep{\partial^2\ue\over\partial t^2},{\partial\over\partial s}(\ue(s)-w_1^\ep(s))\right\rangle_{W',W}\\
&&\qquad\qquad\le c\ep \left\|{\partial\over\partial s}(\ue(s)-w_1^\ep(s))\right\|_W+c\ep^{1/2}\max_{0\le t\le T}\left\|{\partial\over\partial s}(\ue(t)-w_1^\ep(t))\right\|_H\\
&&\qquad\qquad\le c\ep+c\ep^{1/2}\max_{0\le t\le T}\left\|{\partial\over\partial t}(\ue(t)-w_1^\ep(t)\right\|_H\qquad\qquad.
\eeqas
From \eqref{eq:C2} we deduce
\[
\int_0^t\left|\langle\scurl(a^0\curl u_0)-\curl(a^\ep\curl u_1^\ep),{\partial\over\partial s}(\ue(s)-w_1^\ep(s)\rangle_H\right|\le c\ep\left\|{\partial\over\partial s}(\ue-w_1^\ep\right\|_{L^\infty((0,T);W)}\le c\ep
\]
where $c$ is independent of $t$. From \eqref{eq:C5} we get
\beqas
&&\int_0^t\left\langle\scurl(a^\ep\scurl(\ue_1-w^\ep_1)),{\partial\over\partial s}(\ue(s)-w_1^\ep(s))\right\rangle_H ds\\
&&\le c\ep+\\
&& \int_0^t\int_D a^\ep \left((1-\tau^\ep)\curl_yN^r\left(x, \frac{x}{\ep}\right)-\ep \nabla \tau^\ep (x) \times N^r\left(x, \frac{x}{\ep}\right)\right)\cdot(\curl u_0(s,x))_r \frac{\partial  }{\partial s } \curl(u^\ep(s, x)- w_1^\ep(s,x))dxds\\
&&\le c\ep+\\
&&+\int_0^t\int_D a^\ep \left((1-\tau^\ep)\curl_yN^r\left(x, \frac{x}{\ep}\right)-\ep \nabla \tau^\ep (x) \times N^r\left(x, \frac{x}{\ep}\right)\right)\cdot\left[\frac{\partial  }{\partial s }\left(\curl u_0(s,x)_r  \curl(u^\ep(s, x)- w_1^\ep(s,x))\right)\right.\\
&&\left.-\left( \curl(u^\ep(s, x)- w_1^\ep(s,x))\frac{\partial}{\partial s}\curl u_0(s,x)_r\right) \right]dxds.
 \eeqas
Using the initial condition $u_0(0)=w_1^\ep(0)=0$, 
from \eqref{eq:norminDep}, we deduce that
\beqas
&&\left|\int_0^t\left\langle\scurl(a^\ep\scurl(\ue_1-w^\ep_1)),{\partial\over\partial s}(\ue(s)-w_1^\ep(s))\right\rangle_H ds\right|\\
&&\qquad\le c\ep+c\ep^{1/2}\|\curl(\ue-w^\ep_1)(t)\|_H+c\ep^{1/2}\max_{0\le t\le T}\|\curl(\ue-w^\ep_1)(t)\|_H\\
&&\qquad\le c\ep+c\ep^{1/2}\max_{0\le t\le T}\|\curl(\ue-w^\ep_1)(t)\|_H.
\eeqas
Thus
\beqas
&&\int_0^t\left\langle b^\ep \frac{\partial^2 (u^\ep -w_1^\ep)}{\partial s^2}+\curl (a^\ep \curl (u^\ep -w_1^\ep)(s)),{\partial\over\partial s}(\ue-w^\ep_1)(s)\right\rangle ds\\
&&\le c\ep+c\ep^{1/2}\max_{0\le t\le T}\|\curl(\ue-w^\ep_1)(t)\|_H+c\ep^{1/2}\max_{0\le t\le T}\left\|{\partial\over\partial t}(\ue-w^\ep_1)(t)\right\|_H.
\eeqas
The left hand side equals
\beqas
\frac 12\int_0^t{d\over ds}\int_D\left[b^\ep\frac{\partial  (u^\ep- w_1^\ep)}{\partial s }(s)\cdot\frac{\partial  (u^\ep- w_1^\ep)}{\partial s }(s)+a^\ep\curl(\ue(s)-w^\ep_1(s))\cdot(\ue(s)-w^\ep_1(s))\right]dx\\
\ge \alpha\left\|{\partial  (u^\ep- w_1^\ep)\over\partial s }(t)\right\|_H^2-\beta\left\|{\partial  (u^\ep- w_1^\ep)\over\partial s }(0)\right\|_H^2+\alpha\|\curl(\ue(t)-w^\ep_1(t))\|_H
\eeqas
due to condition \eqref{eq:coercive}.

We also have
 \begin{align*}
 \left\|\frac{\partial u^\ep}{\partial t}(0)-\frac{\partial w_1^\ep}{\partial t}(0)\right\|_H&= \left\|\frac{\partial u_0}{\partial t}(0)-\frac{\partial w_1^\ep}{\partial t}(0)\right\|_H \\
 &=\left\|\ep \tau^\ep N^r(\cdot,\frac{\cdot}{\ep})\frac{\partial }{\partial t}\curl u_0(0, \cdot)_r+\ep \nabla \left[\tau^\ep w^r(x,\frac{x}{\ep})\left(\frac{\partial u_{0r}(0,x) }{\partial t}-g_1^r(x)\right)\right]\right\|_H\\ 
 &=\left\|\ep \tau^\ep N^r(x,\frac{x}{\ep})\frac{\partial }{\partial t}\curl u_0(0, \cdot )\right\|_H\\
 &\leq c \ep \left\|\frac{\partial}{\partial t}\curl u_0(0)\right\|_H
 \end{align*}

Thus, we deduce that for all $t \in (0,T)$
 \begin{align*}
 &\|\frac{\partial  (u^\ep- w_1^\ep)}{\partial t }(t, \cdot)\|^2_{H}+\|\curl (u^\ep(t)-w_1^\ep(t))\|^2_H \leq\\
&\qquad c\ep+ 
c\ep^{1/2}\max_{0 \leq t \leq T}\|\frac{\partial  (u^\ep- w_1^\ep)}{\partial t }(t)\|_{H}
 +c\ep^{1/2}\max_{t \in [0,T]}\|\curl (u^\ep(t)-w_1^\ep(t))\|_H
 \end{align*}
 which implies
 \[
\left\|\frac{\partial  (u^\ep- w_1^\ep)}{\partial t }(t)\right\|_{H}+\left\|\curl (u^\ep(t)-w_1^\ep(t))\right\|_H \leq  c\ep^{1/2}.
\]
Due to  ${\partial u_0\over\partial t}\in L^\infty(0,T;H^1(D)^3)$ and $g_1\in H^1(D)^3$, 
\[
\left\|{\partial u_{0r}\over\partial t}\right\|_{L^2(D^\ep)}\le c\ep^{1/2}\ \ \mbox{and}\ \ \|g_{1r}\|_{L^2(D^\ep)}\le c\ep^{1/2}.
\]
 From \eqref{eq:norminDep}
 and the facts that $N^r \in  C^1(\overline{D}, C (\overline{Y})), w^r \in C^1 (\overline{D},C (\overline{Y}) )$ we have 
 \begin{align*}
& \left\|\frac{\partial (u_1^\ep-w_1^\ep)}{\partial t}\right\|_{L^\infty (0, T; H)}\\
&\leq  \left\|\ep (1-\tau^\ep)N^r(x, \frac{x}{\ep}) \frac{\partial}{\partial t} \curl u_{0r}(t)\right.\\
&\qquad\qquad+\left(-\ep \nabla \tau^\ep w^r (x, \frac{x}{\ep})+\ep (1-\tau^\ep)\nabla _xw^r (x, \frac{x}{\ep})
+(1-\tau^\ep)\nabla _yw^r (x, \frac{x}{\ep})\right)\left(\frac{\partial u_{0r}}{\partial t}-g_{1r}\right)\\
&\qquad\qquad\left.+\ep (1-\tau^\ep)w^r (x, \frac{x}{\ep}) \left(\nabla\frac{\partial u_{0r}}{\partial t}-\nabla g_{1r}\right) \right\|_{L^\infty (0, T; H)}\\
&\le c \ep^{1/2}.
 \end{align*}
Similarly, from \eqref{eq:C5}, we have 
\[
\|\curl (u_1^\ep -w_1^\ep)\|_{L^\infty (0,T;H)} \leq c \ep ^{1/2}.
\]
We therefore have
\[
\left\|{\partial(\ue-u_1^\ep)\over\partial t}\right\|_{L^\infty(0,T;H)}+\|\curl(\ue-u_1^\ep)\|_{L^\infty(0,T;H)}\le c\ep^{1/2}.
\]
We note that
\beqas
{\partial u_1^\ep\over\partial t}={\partial u_0\over\partial t}+\ep N^r\left(x,\frac{x}{\ep}\right){\partial\over\partial t}\curl u_0(t,x)+\left(\ep\nabla_xw^r\left(x,{x\over\ep}\right)+\nabla_yw^r\left(x,{x\over\ep}\right)\right)\left({\partial u_{0r}\over\partial t}(t,x)-g_{1r}(x)\right)\\
+\ep w^r\left(x,{x\over\ep}\right)\left(\nabla{\partial u_{0r}\over\partial t}-\nabla g_{1r}\right).
\eeqas
Thus
\[
\left\|{\partial u_1^\ep\over\partial t}-{\partial u_0\over\partial t}-\nabla_yw^r\left(\cdot,{\cdot\over\ep}\right)\left({\partial u_{0r}\over\partial t}(\cdot,\cdot)-g_{1r}(\cdot)\right)\right\|_{L^\infty(0,T;H)}\le c\ep.
\]
Similarly, from \eqref{eq:curlu1ep} we have
\[
\left\|\curl u^\ep_1-\curl u_0-\scurl_yN^r\left(\cdot,{\cdot\over\ep}\right)(\curl u_0)_r\right\|_{L^\infty(0,T;H)}\le c\ep.
\]
We then get the conclusion. \eproof

Next, we derive the homogenization error where $u_0$ only possesses the weaker regularity $W^{2,\infty}(0,T;H^s(\curl,D))$ for $0<s<1$.

\begin{proposition}\label{prop:nonregularhomerror}
Assume that
$g_0=0$, $g_1\in H^1(D)\bigcap W$, $f \in H^1 (0,T; H)$, 
$u_0$ and ${\partial u_0\over\partial t}$ belong to $L^\infty((0,T);H^s(\curl,D))$, and ${\partial^2u_0\over\partial t^2}$ belong to $L^\infty(0,T;H^s(D))$ for $0<s\le 1$, $N^r \in C^1(\overline{D}, C(\overline{Y}))^3$, $\curl_y N^r \in  C^1(\overline{D}, C(\overline{Y}))^3$, $\omega^r \in C^1(\overline{D}, C(\overline{Y}))$ for all $r=1,2,3$. There exists a constant $c$ that does not depend on $\varepsilon$ such that
\[
\left\|\frac{\partial u^\varepsilon}{\partial t}-\left[\frac{\partial u_0}{\partial t}+\nabla_y\frac{\partial \frak u_1}{\partial t}\left(\cdot, \cdot, \frac{\cdot}{\varepsilon} \right ) \right]\right\|_{L^{\infty}(0,T; H)}+\left\|\curl u^\varepsilon-\left[\curl u_0+ \curl_y u_1 \left(\cdot, \cdot, \frac{\cdot}{\varepsilon} \right )\right]\right\|_{L^{\infty}(0,T;H)} \leq c \varepsilon^{\frac{s}{1+s}}.
\]
\end{proposition}

\bproof 
We consider  a set of $M$ open cubes $Q_i~ (i=1,..,M)$ of size $\ep^{s_1}, s_1 >0$ to be chosen later such that $D \subset \bigcup_{i=1}^M Q_i$ and $Q_i \cap D \not= \emptyset$. Each cube $Q_i$ intersects with only a finite number, which does not
depend on $\ep$, of other cubes. We consider a partition of unity that consists of $M$ functions $\rho_i$ such that
$\rho_i$ has support in $Q_i$, $\sum_{i=1}^M \rho_i(x)=1$ for all $x \in D$ and $|\nabla \rho_i| \leq c \ep^{-s_1}$ for all $x$. For $r=1,2,3$ and $i=1,...,M$, we denote by
\[
U^r_i(t)= \frac{1}{|Q_i|} \int_{Q_i} \curl  u_0(t,x)_r dx
\]
and
\[
V^r_i(t)= \frac{1}{|Q_i|} \int_{Q_i}   u_0(t,x)_r dx
\]
(as $u_0 \in H^s(D)^3$ and $\curl u_0 \in H^s(D)^{3}$, for the Lipschitz domain $D$, we can extend each of them,
separately, continuously outside $D$ and understand $u_0$ and curl $u_0$ as these extensions (see Wloka [??]
Theorem 5.6)). Let $U_i$ and $V_i$ denote the vector $(U_i^1, U_i^2,U_i^3)$ and  $(V_i^1, V_i^2,V_i^3)$ respectively. Let $B$ be the unit cube in $\mathbb{R}^3$. From Poincare inequality, we have
$$\int_B|\phi-\int_B \phi(x)dx|^2dx \leq c \int_B|\nabla \phi(x)|^2dx, ~~ \forall \phi \in H^1(B).$$
By translation and scaling, we deduce that
\[
\int_{Q_i}|\phi-\int_{Q_i} \phi(x)dx|^2dx \leq c\ep^{2s_1} \int_{Q_i}|\nabla \phi(x)|^2dx, ~~ \forall \phi \in H^1({Q_i})
\]
i.e
$$\|\phi-\int_{Q_i} \phi(x)dx\|_{L^2(Q_i)} \leq c \ep^{s_1} \|\phi\|_{H^1(Q_i)}.$$
Together with
\[
\|\phi-\int_{Q_i} \phi(x)dx\|_{L^2(Q_i)} \leq c  \|\phi\|_{L^2(Q_i)}
\]
we deduce from interpolation that
\[
\|\phi-\int_{Q_i} \phi(x)dx\|_{L^2(Q_i)} \leq c \ep^{ss_1} \|\phi\|_{H^s(Q_i)}, \forall \phi \in H^s(Q_i).
\]
Thus for $t$ fixed
\begin{equation}
\int_{Q_i}|\curl u_0(t,x)_r-U_i^r(t)|^2 dx \leq c \ep^{2 s s_1}\|(\curl u_0)_r\|^2_{H^s(Q_i)}. 
\label{eq:B}
\end{equation}
We consider the function 
\[
u_1^\ep(t,x)=u_0(t,x)+\ep N^r\left(x,{x\over\ep}\right)U^r_j(t)\rho_j(x)+\ep\nabla\left[ w^r\left(x,{x\over\ep}\right)\left(V^r_j(t)\rho_j(x)-g_1^r(x)t-g_0^r(x)\right)\right ].
\] 
We first show that
\[
\left\|\curl(a^\ep\curl u_1^\ep)-\curl(a^0\curl u_0)\right\|_{L^\infty(0,T;W')} \leq  c (\ep^{1-s_1}+\ep^{ss_1}).
\]
We have
\begin{align}
    \curl u_1^\ep
    =\curl u_0+\ep \curl _xN^r\left(x,{x\over\ep}\right)U^r_j(t)\rho_j(x)+\ep U^r_j(t)\nabla \rho_j(x) \times N^r\left(x,{x\over\ep}\right)+\\
\curl_y N^r\left(x,{x\over\ep}\right)U^r_j(t)\rho_j(x).
\label{eq:1}
\end{align}
Thus
\begin{align*}
a^\ep\curl u_1^\ep(t,x)- a^0\curl u_0(t,x)=&a^\ep\curl u_0+\ep a^\ep\curl _xN^r\left(x,{x\over\ep}\right)U^r_j(t)\rho_j(x)+\ep a^\ep U^r_j(t)\nabla \rho_j(x) \times N^r\left(x,{x\over\ep}\right)\\
&+a^\ep\curl_y N^r\left(x,{x\over\ep}\right)U^r_j(t)\rho_j(x) - a^0(x)\curl u_0(t,x)\\
=&a^\ep U_j(t)\rho_j(x)
+a^\ep\curl_y N^r\left(x,{x\over\ep}\right)U^r_j(t)\rho_j(x) -a^0(x) U_j(t)\rho_j(x)\\ 
&+\ep a^\ep\left[\curl _xN^r\left(x,{x\over\ep}\right)U^r_j(t)\rho_j(x)+ U^r_j(t)\nabla \rho_j(x) \times N^r\left(x,{x\over\ep}\right)\right]\\
&+(a^\ep-a^0)\left(\curl u_0(t,x)- U_j(t)\rho_j(x)\right).
\end{align*}
We have
\begin{multline*}
\curl(a^\ep\curl u_1^\ep)-\curl (a^0\curl u_0)(t,x)=\\
 \curl \left[G_r(x,{x\over\ep})U^r_j(t)\rho_j(x)\right]+\ep\curl I(t,x)+ \curl \left[(a^\ep-a_0)\left(\curl u_0(t,x)- U_j(t)\rho_j(x)\right)\right]
\end{multline*}
where the vector functions $G_r(x,y)$ are defined in \eqref{eq:Gr}
and
\beqas
I(t,x)=a^\ep\left(\curl _xN^r\left(x,{x\over\ep}\right)U^r_j(t)\rho_j(x)+ U^r_j(t)\nabla \rho_j(x) \times N^r\left(x,{x\over\ep}\right)\right).
\eeqas
For all $\phi\in {\cal D}(D)^d$ we have
\begin{multline*}
\langle\curl(a^\ep\curl u_1^\ep)(t)-\curl(a^0\curl u_0)(t),\phi\rangle=\int_DG_r(x,{x\over\ep})U^r_j(t)\rho_j(x)\curl\phi(x)dx\\
+\ep\int_DI(t,x)\curl\phi(x)dx+\int_D(a^\ep-a_0)\left(\curl u_0(t,x)- U^r_j(t)\rho_j(x)\right) \cdot \curl\phi(x)dx\\
=
-\ep\int_D\tilde G_r\left(x,{x\over\ep}\right){\rm div}(U^r_j(t)\rho_j(x)\curl\phi(x))dx-\ep\int_D\nabla_x \tilde G_r\left(x,{x\over\ep}\right)U^r_j(t)\rho_j(x)\curl\phi(x)dx\\+
\ep\int_DI(t,x)\curl\phi(x)dx+\int_D(a^\ep-a_0)\left(\curl u_0(t,x)- U^r_j(t)\rho_j(x)\right) \cdot \curl\phi(x)dx,
\end{multline*}
where $G_r(x,y)=\nabla_y\tilde G_r(x,y)$ as shown above. 
Further,
\[
\int_D\tilde G_r\left(x,{x\over\ep}\right){\rm div}(U^r_j(t)\rho_j(x)\curl\phi(t,x))dx=\int_D\tilde G_r\left(x,{x\over\ep}\right)U^r_j(t)\nabla\rho_j(x)\cdot\curl\phi(t,x).
\]
We note that
$$\left|\int_D \nabla_x \tilde{G}_r \left(x, \frac{x}{\ep}\right)(U^r_j(t)\rho_j) \cdot \curl \phi dx\right| \leq c \|(U^r_j(t)\rho_j)\|_{L^2(D)}\|\curl \phi \|_{L^2(D)^3}.$$
From 
$$\|(U^r_j(t)\rho_j)\|^2_{L^2(D)}= \int_D (U^r_j(t))^2\rho_j(x)^2 dx+ \sum_{i \not=j}\int_D U^r_iU^r_j\rho_i(x)\rho_j(x)dx, $$
and the fact that the support of each function $\rho_i$ intersects only with the support of a finite number
(which does not depend on $\ep$) of other functions $\rho_j$ in the partition of unity, we deduce
\beqas
\|(U^r_j(t)\rho_j)\|^2_{L^2(D)} \leq  c \sum_{j=1}^M (U^r_j(t))^2 |Q_j|= c \sum_{j=1}^M \frac{1}{|Q_j|}\left(\int\limits_{Q_j} \curl u_0(t,x)_r\right)^2 \\
\leq  c \sum_{j=1}^M \int\limits_{Q_j} \curl u_0(t,x)^2_r dx \leq c \int\limits_D\curl u_0(t,x)^2_r dx.
\eeqas
Thus
$$\left|\ep\int_D\nabla_x\tilde G_r\left(x,{x\over\ep}\right)(U^r_j(t)\rho_j)\cdot\curl\phi(t,x))dx\right| \leq c \ep \|\curl \phi \|_{L^2(D)^3}$$
We also have
\begin{align*}
\ep\int\limits_D\tilde G_r\left(x,{x\over\ep}\right){\rm div}\left[U^r_j(t)\rho_j\curl\phi(t,x)\right]dx=&\ep\int\limits_D\tilde G_r\left(x,{x\over\ep}\right)\left[U^r_j(t)\nabla\rho_j(x)\right]\cdot \curl \phi dx\leq&\\
  c \ep \|U^r_j(t)\nabla\rho_j\|_{L^2(D)^3}\|\curl \phi\|_{L^2(D)^3}.
\end{align*}
As the support of each function $\rho_i$ intersects with the support of a finite number of other functions $\rho_j$ and $\|\nabla \rho_j\|_{L^\infty(D)} \leq c \ep^{-s_1}$, we have
\[
\|U^r_j(t)\nabla\rho_j\|^2_{L^2(D)} \leq c \sum_{j=1}^M(U^r_j(t))^2 |Q_j|\|\nabla \rho_j\|^2_{L^\infty(D)} \leq c \ep^{-2s_1}\sum_{j=1}^M(U^r_j(t))^2 |Q_j|\leq c \ep^{-2s_1}
\]
so
\[
\ep\int_D\tilde G_r\left(x,{x\over\ep}\right){\rm div}((U^r_j(t)\rho_j\curl\phi(t,x))dx\leq  c \ep \|U^r_j\nabla\rho_j\|_{L^2(D)^3}\|\curl \phi\|_{L^2(D)^3}\leq c \ep^{1-s_1} \|\curl \phi\|_{L^2(D)^3}.
\]
We have further that
\begin{multline*}
\left \langle\curl \left((a^\ep-a_0)\left(\curl u_0(t,x)- U_j(t)\rho_j(x)\right)\right), \phi  \right \rangle =\int_D(a^\ep-a_0)\left(\curl u_0(t,x)- U_j(t)\rho_j(x)\right) \curl \phi(t,x)dx  \\
\leq c \|\curl u_0(t,x)- U_j(t)\rho_j(x)\|_{L^2(D)^3}\|\curl\phi\|_{L^2(D)^3}.
\end{multline*}
We note that, for $t$ fixed
\[
\int_D |(\curl u_0(t,x))_r- U^r_j(t)\rho_j(x)|^2 dx=\int_D\Big|\sum_{j=1}^M  \left ( (\curl u_0(t,x))_r-U^r_j \right )\rho_j\Big|^2dx.
\]
Using the support property of $\rho_j$, we have from \eqref{eq:B}
\begin{align}
\int_D |(\curl u_0(t,x))_r&- U^r_j(t)\rho_j(x)|^2 dx\leq c\sum_{j=1}^M\int_{Q_j}\Big|  \left ( (\curl u_0(t,x))_r-U^r_j(t) \right )\rho_j\Big|^2dx\nonumber\\
&\leq c \ep^{2 s s_1} \sum_{j=1}^M\|\curl u_0(t,x)_r\|^2_{H^s(Q_j)}\nonumber\\
&=c \ep^{2 s s_1} \sum_{j=1}^M \left [ \int_{Q_j} \curl u_0(t,x)^2_rdx+ \int_{Q_j\times Q_j} \frac{(\curl u_0(t,x)_r-\curl u_0(t,x')_r)^2}{|x-x'|^{2s}}dxdx' \right ]\nonumber\\
&\leq c \ep^{2 s s_1}\left [\|\curl u_0(t,x)_r\|^2_{L^2(D)}+ \int_{D\times D} \frac{(\curl u_0(t,x)_r-\curl u_0(t,x')_r)^2}{|x-x'|^{2s}}dxdx' \right ]\nonumber\\
&=c \ep^{2 s s_1}\|(\curl u_0(t,x))_r\|^2_{H^s(D)}.
\label{eq:curlu0Urj}
\end{align}
Thus
\[
\left \langle\curl \left[(a^\ep-a_0)\left(\curl u_0(t,x)- U_j(t)\rho_j(x)\right)\right], \phi  \right \rangle \leq c \ep^{ss_1}\|\curl\phi\|_{L^2(D)^3}.
\]
From these we  conclude that 
\begin{equation*}
|\langle\curl(a^\ep\curl u_1^\ep)(t)-\curl(a^0\curl u_0)(t),\phi\rangle|\leq c (\ep^{1-s_1}+\ep^{ss_1}) \|\phi\|_W.
\end{equation*}
Using a density argument, we have that this holds for all $\phi\in W$. Thus
\begin{equation}
\|\curl(a^\ep\curl u_1^\ep)-\curl (a^0\curl u_0)\|_{L^\infty (0, T; W')}\leq c (\ep^{1-s_1}+\ep^{ss_1}). \label{D2}
\end{equation}
Since ${\partial u_0\over\partial t}\in L^\infty((0,T);H^s(D))$, by an identical argument, we deduce that
\be
\|\curl(a^\ep\curl {\partial\over\partial t}u_1^\ep)-\curl \left(a^0\curl{\partial\over\partial t} u_0\right)\|_{L^\infty (0, T; W')}\leq c (\ep^{1-s_1}+\ep^{ss_1}).
\label{eq:D2t}
\ee
Choose $s_1=\frac{1}{s+1}$ we have
\begin{equation*}
\|\curl(a^\ep\curl u_1^\ep)-\curl (a^0\curl u_0)\|_{L^\infty (0, T; W')}\le c \ep^{\frac{s}{s+1}}.
\end{equation*}
and
\[
\|\curl(a^\ep\curl {\partial\over\partial t}u_1^\ep)-\curl a^0(\curl{\partial\over\partial t} u_0)\|_{L^\infty (0, T; W')}\leq c\ep^{s\over(s+1)}.
\]
Let $\tau^\ep(x)$ be a function in ${\cal D}(D)$ such that $\tau^\ep(x)=1$ outside an $\ep$ neighbourhood of $\partial D$ and $\sup\limits_{x\in D}\ep|\nabla\tau^\ep(x)|<c$ where $c$ is independent of $\ep$. Let
\begin{equation}
w_1^\ep(t,x)=u_0(t,x)+\ep\tau^\ep(x)N^r\left(x,{x\over\ep}\right)U^r_j(t)\rho_j(x)+\ep\nabla\left[\tau^\ep(x)w^r\left(x,{x\over\ep}\right)\left(V^r_j(t)\rho_j(x)-g_1^r(x)t-g_0^r(x)\right)\right].\label{D3}
\end{equation}
 We then have
\begin{equation}
u_1^\ep-w_1^\ep=\ep(1-\tau^\ep(x))N^r\left(x,{x\over\ep}\right)U^r_j(t)\rho_j(x)+\ep\nabla\left[(1-\tau^\ep(x))w^r\left(x,{x\over\ep}\right)\left(V^r_j(t)\rho_j(x)-g_1^r(x)t-g_0^r(x)\right)\right].\label{D4}
\end{equation}
From this,
\begin{multline}
\curl(u_1^\ep-w_1^\ep)=\ep\curl_xN^r\left(x,{x\over\ep}\right)U^r_j(t)\rho_j(x)(1-\tau^\ep(x))+\curl_yN^r\left(x,{x\over\ep}\right)U^r_j(t)\rho_j(x)(1-\tau^\ep(x))\\
-\ep U^r_j(t)\rho_j(x)\nabla\tau^\ep(x)\times N^r\left(x,{x\over\ep}\right)+\ep(1-\tau^\ep(x))U^r_j(t)\nabla \rho_j(x)\times N^r\left(x,{x\over\ep}\right).\label{D5}
\end{multline}

Let $\tilde{D}^\ep$ be the $3 \ep^{s_1}$ neighbourhood of $\partial D$. We note that $\curl u_0$ is extended continously
 outside $D$. As shown in Hoang and Schwab \cite{HSmultirandom}, for $\phi \in H^1(\tilde{D}^\ep)$
 $$\|\phi\|_{L^2(\tilde{D}^\ep)} \leq c \ep^{s_1/2}\|\phi\|_{H^1(\tilde{D}^\ep)}.$$
 From this and 
 $$\|\phi\|_{L^2(\tilde{D}^\ep)} \leq \|\phi\|_{L^2(\tilde{D}^\ep)} ,$$
 using interpolation we get
 \begin{equation}
 \|\phi\|_{L^2(\tilde{D}^\ep)} \leq c \ep^{ss_1/2}\|\phi\|_{H^s(\tilde{D}^\ep)} \leq c \ep^{ss_1/2}\|\phi\|_{H^s(D)}, \label{D7}
 \end{equation}
 for all $\phi \in H^s(D)$ extended continuously outside $D$.
 We then have for $t$ fixed,
 \begin{align*}
     \|U^r_j(t) \rho_j\|^2_{L^2(D^\ep)} &\leq c \sum_{j=1}^M \int_{Q_j \cap D^\ep} (U^r_j(t))^2 \rho_j^2 dx\\
     &\leq c \sum_{j=1}^M |Q_j \cap D^\ep| \frac{1}{|Q_j|^2} \left(\int_{Q_j}(\curl u_0)_r dx\right)^2\\
     &\leq c \sum_{Q_j \cap D^\ep \not= \emptyset}  \frac{|Q_j \cap D^\ep|}{|Q_j|} \int_{Q_j}(\curl u_0)^2_r dx.
 \end{align*}
 As $\partial D$ is Lipschitz,  $D^\ep$ is the $\ep$ neighbourhood of $\partial D$ and $Q_j$ has size $\ep^{s_1}$, $|Q_j \cap D^\ep|  \leq c \ep^{1+(d-1)s_1}$ so $ \frac{|Q_j \cap D^\ep|}{|Q_j|} \leq c \ep^{1-s_1}$.  When $Q_j \cap D^\ep \not= \emptyset$ , $Q_j \subset \tilde{D}^\ep.$ Thus 
 \be
\|U^r_j(t) \rho_j\|^2_{L^2(D^\ep)} \leq c \ep^{1-s_1} \|(\curl u_0(t,x))_r\|^2_{L^2(\tilde{D}^\ep)} \leq c \ep^{1-s_1+ss_1} \|\curl u_0(t,x)\|^2_{H^s(D)^3}
\label{eq:es1}
\ee
 Therefore
\[
\|\curl_yN^r\left(x,{x\over\ep}\right)U^r_j(t)\rho_j(x)(1-\tau^\ep(x))\|_{L^2(D^\ep)^3} \leq c \ep^{(1-s_1+ss_1)/2}
\]
 and 
\[
\|\ep U^r_j(t)\rho_j(x)\nabla\tau^\ep(x)\times N^r\left(x,{x\over\ep}\right)\|_{L^2(D^\ep)^3} \leq c \ep^{(1-s_1+ss_1)/2}.
\] 
Similarly, we have
 \begin{align}
     \|U^r_j(t) \nabla \rho_j\|^2_{L^2(D^\ep)^3} &\leq c \ep^{-2s_1} \sum_{j=1}^M |Q_j \cap D^\ep| (U^r_j(t))^2\nonumber\\
     &\leq c \ep^{-2s_1}\sum_{Q_j \cap D^\ep \not= \emptyset}  \frac{|Q_j \cap D^\ep|}{|Q_j|} \int_{Q_j}(\curl u_0)^2_r dx\nonumber\\
     &\leq c \ep^{-2s_1+1-s_1}\|\curl u_0\|^2_{L^2(\tilde{D}^\ep)}\nonumber\\
     &\leq c \ep^{1-3s_1+ss_1}\|\curl u_0\|^2_{H^s(D)^3}.\label{eq:es2}
 \end{align}
Thus
\[
\|\ep(1-\tau^\ep(x))U^r_j(t)\nabla \rho_j(x)\times N^r\left(x,{x\over\ep}\right)\|_{L^2(D)^3} \leq c \ep^{(1-s_1)+(1-s_1+ss_1)/2}.
\]
Therefore
\be
\|\curl(u_1^\ep-w_1^\ep)\|_{L^2(D)^3} \leq c(\ep^{(1-s_1+ss_1)/2}+ \ep^{(1-s_1)+(1-s_1+ss_1)/2}).
\label{eq:u1epw1ep}
\ee
Arguing as above   (with $t$ fixed) we deduce that
\[
\|V^r_j(t)\rho_j\|_{L^2(D^\ep)} \leq c \ep^{(1-s_1+ss_1)/2}, ~\|V^r_j(t)\nabla\rho_j\|_{L^2(D^\ep)} \leq c \ep^{(1-s_1+ss_1)/2-s_1} .
\]
 We have
 \begin{align*}
     \frac{\partial^2 w_1^\ep}{\partial  t^2}=&\frac{\partial^2 u_0}{\partial  t^2}+\ep \tau^\ep N^r\left(x,{x\over\ep}\right) \frac{\partial ^2U^r_j}{\partial t^2}(t)\rho_j(x)+ \ep \nabla \tau^\ep w^r\left(x,{x\over\ep}\right) \frac{\partial ^2  V^r_j}{\partial t^2}(t)\rho_j(x)\\
     &+\ep\tau^\ep \nabla_x  w^r\left(x,{x\over\ep}\right) \frac{\partial ^2  V^r_j}{\partial t^2}(t)\rho_j(x)+\tau^\ep \nabla_y w^r\left(x,{x\over\ep}\right) \frac{\partial ^2  V^r_j}{\partial t^2}(t)\rho_j(x)+\ep \tau^\ep w^r\left(x,{x\over\ep}\right) \frac{\partial ^2V^r_j  }{\partial t^2}(t)\nabla\rho_j(x).
 \end{align*}
 Thus
 \begin{align*}
    b^0 \frac{\partial^2 u_0 }{\partial t^2}-& b^\ep\frac{\partial^2 w_1^\ep}{\partial  t^2}= b^0 \frac{\partial^2 u_0 }{\partial t^2}-b^\ep\frac{\partial^2 u_0}{\partial  t^2}-b^\ep\ep \tau^\ep N^r\left(x,{x\over\ep}\right) \frac{\partial ^2 U^r_j}{\partial t^2}(t)\rho_j(x)- b^\ep\ep \nabla \tau^\ep w^r\left(x,{x\over\ep}\right) \frac{\partial ^2  V^r_j}{\partial t^2}(t)\rho_j(x)\\
     &-b^\ep\ep\tau^\ep \nabla_x  w^r\left(x,{x\over\ep}\right) \frac{\partial ^2  V^r_j}{\partial t^2}(t)\rho_j(x)-b^\ep\tau^\ep \nabla_y w^r\left(x,{x\over\ep}\right) \frac{\partial ^2  V^r_j}{\partial t^2}(t)\rho_j(x)-b^\ep\ep \tau^\ep w^r\left(x,{x\over\ep}\right) \frac{\partial ^2V_j^r  }{\partial t^2}(t) \nabla \rho_j(x)\\
     &=b^0\frac{\partial ^2  V_j}{\partial t^2}(t)\rho_j(x)-b^\ep\frac{\partial ^2  V_j}{\partial t^2}(t)\rho_j(x)-b^\ep \nabla_y w^r\left(x,{x\over\ep}\right) \frac{\partial ^2  V^r_j}{\partial t^2}(t)\rho_j(x)\\
     &-b^\ep\ep \tau^\ep N^r\left(x,{x\over\ep}\right) \frac{\partial ^2 U^r_j}{\partial t^2}(t)\rho_j(x)- b^\ep\ep \nabla \tau^\ep w^r(x, \frac{x}{\ep}) \frac{\partial ^2  V^r_j}{\partial t^2}(t)\rho_j(x)-b^\ep\ep\tau^\ep \nabla_x  w^r\left(x,{x\over\ep}\right) \frac{\partial ^2  V^r_j}{\partial t^2}(t)\rho_j(x)\\
     &+b^\ep(1-\tau^\ep) \nabla_y w^r\left(x,{x\over\ep}\right) \frac{\partial ^2  V^r_j}{\partial t^2}(t)\rho_j(x)-b^\ep\ep \tau^\ep w^r\left(x,{x\over\ep}\right) \frac{\partial ^2 V_j^r }{\partial t^2}(t)\nabla \rho_j(x)\\
     &+ (b^0-b^\ep)\left( \frac{\partial^2 u_0 }{\partial t^2}(t)-\frac{\partial ^2 V_j}{\partial t^2}(t)\rho_j(x)  \right).
\end{align*}
For $\phi\in W$, we consider
 \begin{align*}
    J:&= \left\langle b^0\frac{\partial ^2  V_j}{\partial t^2}(t)\rho_j(x)-b^\ep\frac{\partial ^2  V_j}{\partial t^2}(t)\rho_j(x)-b^\ep(x) \nabla_y w^r\left(x,{x\over\ep}\right) \frac{\partial ^2  V^r_j}{\partial t^2}(t)\rho_j(x), \phi\right\rangle_{W',W}\\
     &=\int_D b^0_{ir}\frac{\partial ^2  V^r_j}{\partial t^2}(t)\rho_j(x)\phi_i-b^\ep_{ik}(x)\left(\delta_{kr}+\frac{\partial w^r}{\partial y_k}\left(x,{x\over\ep}\right)\right)\frac{\partial ^2  V^r_j}{\partial t^2}(t)\rho_j(x)\phi_idx\\
&=\int_D F_r(x,{x\over\ep}){\partial^2V_j^r\over\partial t^2}(t)\rho_j(x)dx
 \end{align*}
where the function $F_r(x,y)$ is defined in \eqref{eq:F}. 
Using the function $\tilde F_r$ defined in \eqref{eq:F1}, we have 
 \begin{align*}
J&=\int_D \Big[\ep \curl \tilde F_r(x,\frac{x}{\ep})-\ep \curl_x\tilde F_r(x,\frac{x}{\ep})\Big]\frac{\partial ^2  V^r_j}{\partial t^2}(t)\rho_j(x)\cdot \phi dx.\\
&=\ep \int_D \tilde F_r(x,\frac{x}{\ep})\cdot \curl \big(\frac{\partial ^2  V^r_j}{\partial t^2}(t)\rho_j(x) \phi \big) dx-\ep \int_D \curl_x\tilde F_r(x,\frac{x}{\ep})\frac{\partial ^2  V^r_j}{\partial t^2}(t)\rho_j(x) \cdot \phi dx\\
&\le\ep \int_D \tilde F_r(x,\frac{x}{\ep})\left(\curl \phi \frac{\partial ^2  V^r_j}{\partial t^2}(t)\rho_j(x)+ \phi \times \frac{\partial ^2  V^r_j}{\partial t^2}(t)\nabla \rho_j(x)\right)dx+ c\ep\|\phi\|_{H}\\
&\le c\ep\|\curl \phi\|_{H}+ c\ep^{1-s_1} \|\phi\|_{H} + c\ep\|\phi\|_{H}
 \end{align*}
As ${\partial^2u_0\over\partial t^2}\in L^\infty(0,T;H^s(D))$, there is a constant $c$ such that for all $t\in (0,T)$
\beqas
 \left |\int_D(b^\ep-b^0)\left( \frac{\partial^2 u_0 }{\partial t^2}-\frac{\partial ^2 V_j(t) }{\partial t^2}\rho_j(x) \right) \cdot \phi dx \right | &\leq& c \left\|\sum_{j=1}^M (\frac{\partial^2 u_0 }{\partial t^2}-\frac{\partial ^2 V_j(t) }{\partial t^2})\rho_j\right\|_{H}\| \phi\|_{H}\\
& \leq&  c \ep^{ss_1} \|\phi\|_{H}.
\eeqas
Thus 
\begin{eqnarray}
&&\left\langle b^0 \frac{\partial^2 u_0 }{\partial t^2}- b^\ep\frac{\partial^2 w_1^\ep}{\partial  t^2} , \phi\right\rangle_{W',W}\\
&&\qquad\qquad\le c\ep\|\curl \phi\|_{(L^2(D))^d}+ c\ep^{1-s_1} \|\phi\|_{(L^2(D))^d}+
c\ep^{\frac{1-s_1+ss_1}{2}}\|\phi\|_{(L^2(D))^d}+c\ep^{ss_1}\|\phi\|_{(L^2(D))^d}\nonumber\\
&&\qquad\qquad\le c\ep\|\curl\phi\|_{L^2(D)^d}+c\ep^{s\over 1+s}\|\phi\|_{L^2(D)^3}
\label{eq:bdist}
\end{eqnarray}
when we choose $s_1=1/(1+s)$. 
Using 
\[ 
b^\ep \frac{\partial u^\ep }{\partial t^2}+\curl (a^\ep \curl u^\ep)=b^0 \frac{\partial^2 u_0 }{\partial t^2}+\curl (a^0 \curl u_0)
\]
we have
 \begin{align*}
     b^\ep \frac{\partial (u^\ep -w_1^\ep)}{\partial t^2}&+\curl (a^\ep \curl (u^\ep -w_1^\ep))\\
     &=b^0 \frac{\partial^2 u_0 }{\partial t^2}- b^\ep \frac{\partial^2 w_1^\ep }{\partial t^2}+\curl (a^\ep\curl(u_1^\ep-w_1^\ep))+\curl (a^0\curl u_0)-\curl (a^\ep\curl u_1^\ep).
\end{align*}
As shown in the proof of Proposition \ref{prop:2sregularhomerror}, ${\partial\ue\over\partial t}$ is uniformly bounded in $L^\infty(0,T;W)$ with respect to $\ep$. For $w_1^\ep$, we have
\beqas
&&{\partial w^\ep_1\over\partial t}={\partial u_0\over\partial t}(t,x)+\ep\tau^\ep(x)N^r\left(x,{x\over\ep}\right){\partial U_j^r\over\partial t}(t)\rho_j(x)\\
&&\qquad+\left(\ep\nabla\tau^\ep(x) w^r\left(x,{x\over\ep}\right)+\ep\tau^\ep(x)\nabla_xw^r\left(x,{x\over\ep}\right)+\tau^\ep(x)\nabla_y w^r\left(x,{x\over\ep}\right)\right)\left({\partial V_j^r\over\partial t}(t)\rho_j(x)-g_{1r}\right)\\
&&\qquad+\ep\tau^\ep(x)w^r\left(x,{x\over\ep}\right)\left({\partial V_j^r\over\partial t}(t)\nabla\rho_j(x)-\nabla g_{1r}(x)\right)
\eeqas
and 
\beqas
&&\curl{\partial w^\ep_1\over\partial t}=\curl{\partial u_0\over\partial t}(t,x)+\tau^\ep(x){\partial U_j^r\over\partial t}(t)\rho_j(x)\left(\ep\scurl_xN^r\left(x,{x\over\ep}\right)+\scurl_yN^r\left(x,{x\over\ep}\right)\right)\\
&&\qquad\qquad\qquad+{\partial U_j^r\over\partial t}(t)\left(\ep\nabla\tau^\ep(x)\rho_j(x)+\ep\tau^\ep(x)\nabla\rho_j(x)\right)\times N^r\left(x,{x\over\ep}\right). 
\eeqas
As ${\partial u_0\over\partial t}\in L^\infty(0,T;H(\curl,D))$, $\|{\partial U_j^r\over\partial t}\rho_j\|_{L^2(D)^3}\le c$, $\|{\partial V_j^r\over\partial t}\rho_j\|_{L^2(D)^3}\le c$, $\|{\partial U_j^r\over\partial t}\nabla\rho_j\|_{L^2(D)^3}\le c\ep^{-s_1}$ and $\|{\partial V_j^r\over\partial t}\nabla\rho_j\|_{L^2(D)^3}\le c\ep^{-s_1}$. Therefore $\curl {\partial w^\ep_1\over\partial t}$ and ${\partial w^\ep_1\over\partial t}$ are uniformly bounded in $H$ with respect to $\ep$, i.e. ${\partial w^\ep_1\over\partial t}$ is uniformly bounded in $L^\infty(0,T;H^1(\curl,D))$.
We have
\beqas
&&\int_0^t\left\langle b^0 \frac{\partial^2 u_0 }{\partial t^2}(s)- b^\ep \frac{\partial^2 w_1^\ep }{\partial t^2}(s),{\partial(u^\ep-w_1^\ep)\over\partial s}(s)\right\rangle_H\\
&&\qquad\le \int_0^t\left(c\ep\left\|\curl{\partial(u^\ep-w_1^\ep)\over\partial s}(s)\right\|_{H}+c\ep^{s\over 1+s}\sup_{0\le s\le T}\left\|{\partial(u^\ep-w_1^\ep)\over\partial s}(s)\right\|_{H}\right)dt\\
&&\qquad\le c\ep+c\ep^{s\over 1+s}\sup_{0\le t\le T}\left\|{\partial(u^\ep-w_1^\ep)\over\partial s}(t)\right\|_{H}.
\eeqas
We also have
\beqas
&&\int_0^t\left\langle\curl (a^0\curl u_0(s))-\curl (a^\ep\curl u_1^\ep(s)),{\partial(u^\ep-w_1^\ep)\over\partial s}(s)\right\rangle_{W',W}ds=\\
&&\qquad\int_0^t{\partial\over\partial s}\left\langle\curl (a^0\curl u_0(s))-\curl (a^\ep\curl u_1^\ep(s)),(u^\ep-w_1^\ep)(s)\right\rangle_{W',W}ds-\\
&&\qquad\int_0^t\left\langle\curl (a^0\curl {\partial u_0\over\partial s}(s))-\curl (a^\ep\curl {\partial u_1^\ep\over\partial s}(s)),(u^\ep-w_1^\ep)(s)\right\rangle_{W',W}ds.
\eeqas
Since $u_0(0)=0$ and $u_1^\ep(0)=0$, together with \eqref{eq:D2t} we have that
\beqas
&&\left|\int_0^t\left\langle\curl (a^0\curl u_0(s))-\curl (a^\ep\curl u_1^\ep(s)),{\partial(u^\ep-w_1^\ep)\over\partial s}(s)\right\rangle_H\right|ds\\
&& \qquad\le\left|\left\langle \curl(a^0(\curl u_0(t))-\curl(a^\ep\curl u_1^\ep(t)),(\ue-w^\ep_1)(t)\right\rangle\right|+c\int_0^t\ep^{s\over s+1}\|(\ue-w^\ep_1)(s)\|_V\\
&&\qquad\le c\ep^{s\over s+1}\|(\ue-w^\ep_1)(t)\|_V+c\ep^{s\over s+1}\sup_{0\le t\le T}\|(\ue-w^\ep_1)(t)\|_V.
\eeqas
Now we estimate
\[
\int_0^t\left\langle\curl(a^\ep\curl(u_1^\ep(s)-w^\ep_1(s)),{\partial(\ue-w^\ep_1)\over\partial s}(s)\right\rangle_Hds
\]
using \eqref{D5}. We have that
\beqas
&&\left|\int_0^t\left\langle\ep\curl\left(a^\ep\scurl_xN^r\left(x,{x\over\ep}\right)U^r_j(t)\rho_j(x)(1-\tau^\ep(x))\right),{\partial(\ue-w^\ep_1)\over\partial s}(s)\right\rangle_{W',W} ds\right|\\
&&\qquad\qquad\le c\ep\int_0^t\left\|\curl{\partial(\ue-w^\ep_1)\over\partial s}(s)\right\|_Vds\le c\ep;
\eeqas
and
\beqas
&&\left|\int_0^t\left\langle\ep\curl\left(a^\ep(1-\tau^\ep(s))U_j^r(t)\nabla\rho_j(x)\times N^r\left(x,{x\over\ep}\right)\right),{\partial(\ue-w^\ep_1)\over\partial s}(s)\right\rangle_{W',W} dx\right|\\
&&\qquad\qquad\le c\ep\int_0^t\|U_j^r(s)\nabla\rho_j\|_{L^2(D^\ep)^3}\left\|{\partial(\ue-w^\ep_1)\over\partial s}(s)\right\|_V\le c\ep\ep^{{1-3s_1+ss_1\over 2}}=c\ep^{2s\over s+1}.
\eeqas
For the other two terms in \eqref{D5}, we have
\beqas
&&\int_0^t\left\langle\curl\left(a^\ep\scurl_yN^r\left(x,{x\over\ep}\right)U_j^r(s)\rho_j(x)(1-\tau^\ep(x))\right),{\partial(\ue-w^\ep_1)\over\partial s}(s)\right\rangle_{W',W}ds\\
&&=\int_0^t{\partial\over\partial s}\int_Da^\ep(x)\scurl_yN^r\left(x,{x\over\ep}\right)U_j^r(s)\rho_j(x)(1-\tau^\ep(x))\cdot\curl(\ue(s)-w^\ep_1(s))dxds-\\
&&\qquad\int_0^t\int_Da^\ep(x)\scurl_yN^r\left(x,{x\over\ep}\right){\partial U_j^r\over\partial s}(s)\rho_j(x)(1-\tau^\ep(x))\cdot\curl(\ue(s)-w^\ep_1(s))ds\\
&&=\int_Da^\ep\scurl_yN^r\left(x,{x\over\ep}\right)U_j^r(t)\rho_j(x)(1-\tau^\ep(x))\curl(\ue(t)-w^\ep_1(t))-\\
&&\qquad\int_0^t\int_Da^\ep(x)\scurl_yN^r\left(x,{x\over\ep}\right){\partial U_j^r\over\partial s}(s)\rho_j(x)(1-\tau^\ep(x))\cdot\curl(\ue(s)-w^\ep_1(s))ds.
\eeqas
Thus
\beqas
&&\left|\int_0^t\left\langle\curl\left(a^\ep\scurl_yN^r\left(x,{x\over\ep}\right)U_j^r(s)\rho_j(x)(1-\tau^\ep(x))\right),{\partial(\ue-w^\ep_1)\over\partial s}(s)\right\rangle_{W',W}ds\right|\\
&&\le\|U_j^r(t)\rho_j\|_{L^2(D^\ep)}\|\curl(\ue(t)-w^\ep_1(t))\|_{L^2(D)^3}+c\int_0^t\left\|{\partial U^r_j(s)\over\partial s}\rho_j\right\|_{L^2(D^\ep)}\|\curl(\ue(s)-w^\ep_1(s))\|_{L^2(D)^3}ds\\
&&\le c\ep^{1-s_1+ss_1\over 2}\sup_{0\le t\le T}\|\curl(\ue(t)-w^\ep_1(t))\|_{L^2(D)^3}
\eeqas
as $u_0$ and ${\partial u_0\over\partial t}$ belong to $L^\infty(0,T;H^s(D)))$. 
Similarly, we have
\beqas
&&\left|\int_0^t\left\langle\ep\curl\left(a^\ep U_j^r(s)\rho_j(x)\nabla\tau^\ep(x)\times N^r\left(x,{x\over\ep}\right)\right),{\partial(\ue-w^\ep_1)\over\partial s}(s)\right\rangle ds\right|\\
&&\qquad\qquad\le c\ep^{1-s_1+ss_1\over 2}\sup_{0\le t\le T}\|\curl(\ue(t)-w^\ep_1(t))\|_{L^2(D)^3}.
\eeqas
Therefore,
\[
\left|\int_0^t\left\langle\curl(a^\ep\curl(u^\ep_1(s)-w^\ep_1(s)),{\partial(\ue-w^\ep_1)\over\partial s}(s)\right\rangle_H\right|\le c\ep^{2s\over s+1}+c\ep^{s\over s+1}\sup_{0\le t\le T}\|\curl(\ue(t)-w^\ep_1(t))\|_{L^2(D)^3}.
\]
We then deduce
\beqas
&&\int_0^t\left\langle b^\ep{\partial^2(\ue-w^\ep_1)\over\partial s^2}(s)+\curl(a^\ep\curl(\ue-w^\ep_1))(s),{\partial(\ue-w^\ep_1)\over\partial s}(s)\right\rangle\\
&&\qquad\qquad\le c\ep+c\ep^{2s\over s+1}+c\ep^{s\over 1+s}\sup_{0\le t\le T}\left\|{\partial(\ue-w^\ep_1)\over\partial t}(t)\right\|_{L^2(D)^3}+c\ep^{s\over 1+s}\sup_{0\le t\le T}\left\|\curl(\ue(t)-w^\ep_1(t)\right\|_{L^2(D)^3}.
\eeqas
Therefore
\beqas
&&\frac12\int_Db^\ep(x){\partial(\ue-w^\ep_1)\over\partial t}(t)\cdot{\partial(\ue-w^\ep_1)\over\partial t}(t)dx+\frac12\int_Da^\ep(x)\curl(\ue-w^\ep_1)(t)\cdot\curl(\ue-w^\ep_1)(t)dx\\
&&\qquad\le c\ep^{2s\over s+1}+c\ep^{s\over 1+s}\sup_{0\le t\le T}\left\|{\partial(\ue-w^\ep_1)\over\partial t}(t)\right\|_{L^2(D)^3}+c\ep^{s\over 1+s}\sup_{0\le t\le T}\left\|\curl(\ue(t)-w^\ep_1(t))\right\|_{L^2(D)^3}\\
&&\qquad+\int_Db^\ep(x){\partial(\ue-w^\ep_1)\over\partial t}(0)\cdot{\partial(\ue-w^\ep_1)\over\partial t}(0)dx
\eeqas
(note that $\ue(0)=w^\ep_1(0)=0$). We have 
 \begin{align*}
  &\left\|\frac{\partial u^\ep}{\partial t}(0)-\frac{\partial w_1^\ep}{\partial t}(0)\right\|_H\\
&\qquad=  \left\|\frac{\partial u_0}{\partial t}(0)-\frac{\partial w_1^\ep}{\partial t}(0)\right\|_H  \\
 &\qquad=\left\|\ep \tau^\ep N^r\left(x,\frac{x}{\ep}\right)\frac{\partial U^r_j}{\partial t}(0)\rho_j(x)+\ep \nabla\left(\tau^\ep w^r\left(x,\frac{x}{\ep}\right)\left(\frac{\partial V^r_j}{\partial t}(0)\rho_j(x)-g_1^r(x)\right)\right)\right\|_H\\
 &\qquad = \left\|\ep \tau^\ep N^r\left(x,\frac{x}{\ep}\right)\frac{\partial U^r_j}{\partial t}(0)\rho_j(x)\right.\\
&\qquad\qquad+\left(\ep \nabla \tau^\ep w^r\left(x,\frac{x}{\ep}\right)+\ep \tau^\ep \nabla _xw^r\left(x,\frac{x}{\ep}\right)
  + \tau^\ep \nabla _yw^r\left(x,\frac{x}{\ep}\right)\right)\left(\frac{\partial V^r_j}{\partial t}(0)\rho_j(x)-g_{1r}(x)\right)\\
&\qquad\qquad   +\left.\ep \tau^\ep  w^r\left(x,\frac{x}{\ep}\right)\left(\frac{\partial V^r_j}{\partial t}(0)\nabla\rho_j(x)-\nabla g_{1r}(x)\right)\right\|_H\\
&\le c\ep+c\left\|\frac{\partial V^r_j}{\partial t}(0)\rho_j(x)-g_{1r}(x)\right\|_{H}+
c\ep\left\|{\partial V_j^r\over\partial t}(0)\nabla\rho_j-\nabla g_{1r}\right\|_H.
\end{align*}
As $g_1\in H^1(D)^3$, a similar argument as for showing \eqref{eq:curlu0Urj} for $s=1$ shows that
\[
\left\|\frac{\partial V^r_j}{\partial t}(0)\rho_j(x)-g_{1r}(x)\right\|_{H}\le\ep^{s_1}<c\ep^{s\over 1+s}.
\]
Further,
\[
\left\|\frac{\partial V^r_j}{\partial t}(0)\nabla\rho_j(x)\right\|_{L^2(D^\ep)^3}\le c\ep^{-s_1}.
\]
Thus
\beqas
\left\|\frac{\partial u^\ep}{\partial t}(0)-\frac{\partial w_1^\ep}{\partial t}(0)\right\|_H\leq c\ep^{s\over 1+s}+c\ep^{1-s_1}\le c\ep^{s\over s+1}.
 \eeqas
Using \eqref{eq:coercive} we get
\beqas
&&\left\|{\partial(\ue-w^\ep_1)\over\partial t}(t)\right\|_H^2+\left\|\curl(\ue(t)-w^\ep_1(t))\right\|_H^2\\
&&\qquad\qquad\le c\ep^{2s\over s+1}+c\ep^{s\over s+1}\max_{0\le t\le T}\left\|\frac{\partial(u^\ep-w_1^\ep)}{\partial t}\right\|_H+c\ep^{s\over s+1}\max_{0\le t\le T}\|\curl{u^\ep(t)-w_1^\ep(t)}\|_H.
\eeqas
From this we deduce that for all $t\in (0,T)$
\be
\left\|{\partial(\ue-w^\ep_1)\over\partial t}(t)\right\|_H+\left\|\curl(\ue(t)-w^\ep_1(t))\right\|_H\le c\ep^{s\over s+1}.
\label{eq:uew1ep}
\ee
From \eqref{D4} we have
\beqas
&&{\partial(u_1^\ep-w^\ep_1)\over\partial t}(t)=\ep(1-\tau^\ep(x))N^r\left(x,{x\over\ep}\right){\partial U_j^r\over\partial t}(t)\rho_j(x)\\
&&+\left[-\ep\nabla\tau^\ep(x)w^r\left(x,{x\over\ep}\right)+\ep(1-\tau^\ep(x))\nabla_xw^r\left(x,{x\over\ep}\right)+(1-\tau^\ep(x))\nabla_yw^r\left(x,{x\over\ep}\right)\right]\left({\partial V_j^r\over\partial t}(t)\rho_j(x)-g_{1r}\right)\\
&&+\ep(1-\tau^\ep(x))w^r\left(x,{x\over\ep}\right)({\partial V_j^r\over\partial t}(t)\nabla \rho_j(x)-\nabla g_{1r}).
\eeqas
Therefore, using $g_1\in H^1(D)^3$ we get
\beqas
\left\|{\partial(u_1^\ep-w^\ep_1)\over\partial t}(t)\right\|_H\le c\ep\left\|{\partial U_j^r\over\partial t}(t)\rho_j\right\|_{L^2(D^\ep)^3}+\left\|{\partial V_j^r\over\partial t}(t)\rho_j\right\|_{L^2(D^\ep)^3}+c\|g_1\|_{L^2(D^\ep)^3}+c\ep\left\|{\partial V_j^r\over\partial t}(t)\nabla\rho_j-\nabla g_1^r\right\|_{L^2(D^\ep)^3}.
\eeqas
As ${\partial\over\partial t}\curl u_0\in L^\infty(0,T;H^s(D)^3)$, ${\partial\over\partial t}u_0\in L^\infty(0,T;H^s(D)^3)$ and $g_1\in H^1(D)^3$, we deduce that 
\beqas
\left\|{\partial U_j^r\over\partial t}(t)\rho_j\right\|_{L^2(D^\ep)^3}\le c,\ \ \left\|{\partial V_j^r\over\partial t}(t)\rho_j\right\|_{L^2(D^\ep)^3}\le c\ep^{1-s_1+ss_1\over 2},\\
 \left\|{\partial V_j^r\over\partial t}(t)\nabla\rho_j\right\|_{L^2(D^\ep)^3}\le c\ep^{1-3s_1+ss_1\over 2},\ \ \|g_1\|_{L^2(D^\ep)^3}\le c\ep^{1/2}.
\eeqas
Therefore
\be
\left\|{\partial(u_1^\ep-w^\ep_1)\over\partial t}(t)\right\|_H\le c\ep^{1-s_1+ss_1\over 2}+\ep\ep^{1-3s_1+ss_1\over 2}+c\ep^{1/2}\le c\ep^{s\over s+1}.
\label{eq:dtu1epw1ep}
\ee
From \eqref{eq:u1epw1ep}, \eqref{eq:uew1ep} and \eqref{eq:dtu1epw1ep}, we deduce that
\be
\left\|{\partial(\ue-u^\ep_1)\over\partial t}(t)\right\|_H+\left\|\curl(\ue(t)-u^\ep_1(t))\right\|_H\le c\ep^{s\over s+1}.
\label{eq:i1}
\ee
We note that
\[
\left\|\ep\scurl_xN^r\left(x,{x\over\ep}\right)U^r_j\rho_j\right\|_H\le c\ep,\ \mbox{ and }\ \left\|\ep N^r\left(x,{x\over\ep}\right)\times (U^r_j\nabla\rho_j)\right\|_H\le c\ep\ep^{-s_1}=c\ep^{s\over s+1}
\]
so from \eqref{eq:1}
\[
\|\curl u_1^\ep-[\curl u_0+\scurl_yN^r\left(x,{x\over\ep}\right)U_j^r\rho_j]\|_H\le c\ep^{s\over s+1}.
\]
Thus we deduce from \eqref{eq:curlu0Urj} that
\be
\|\curl u_1^\ep-[\curl u_0+\scurl_yN^r\left(x,{x\over\ep}\right)\curl u_0(x)_r]\|_H\le c\ep^{s\over s+1}.
\label{eq:i2}
\ee
We further have
\beqas
{\partial u_1^\ep\over\partial t}(t)={\partial u_0\over\partial t}(t)+\ep N^r\left(x,{x\over\ep}\right){\partial U_j^r\over\partial t}(t)\rho_j(x)+
\left(\ep\nabla_x w^r\left(x,{x\over\ep}\right)+\nabla_yw^r\left(x,{x\over\ep}\right)\right)\left({\partial V_j^r\over\partial t}\rho_j(x)-g_{1r}\right)+\\
\ep w^r\left(x,{x\over\ep}\right)\left({\partial V_j^r(t)\over\partial t}\nabla\rho_j(x)-\nabla g_{1r}(x)\right).
\eeqas
As 
\beqas
\left\|N^r\left(\cdot,{\cdot\over\ep}\right){\partial U_j^r\over\partial t}(t)\rho_j\right\|_H\le c,\ \ \left\|{\partial V_j^r\over\partial t}\rho_j\right\|_H\le c,\ \ \left\|{\partial V_j^r(t)\over\partial t}\nabla\rho_j\right\|_H\le c\ep^{-s_1}
\eeqas
we have that
\beqas
\left\|{\partial u_1^\ep\over\partial t}(t)-{\partial u_0\over\partial t}(t)-\nabla_yw^r\left(x,{x\over\ep}\right)\left({\partial V_j^r\over\partial t}(t)\rho_j-g_{1r}\right)\right\|_H\le c\ep^{s\over s+1}.
\eeqas
Using
\[
\left\|{\partial V^r_j\over\partial t}(t)\rho_j(x)-{\partial u_{0r}\over\partial t}(t)\right\|_H\le c\ep^{ss_1}=c\ep^{s\over s+1}
\]
we deduce that
\be
\left\|{\partial u_1^\ep\over\partial t}(t)-{\partial u_0\over\partial t}(t)-\nabla_yw^r\left(x,{x\over\ep}\right)\left({\partial u_{0r}\over\partial t}(t)-g_{1r}(x)\right)\right\|_H\le c\ep^{s\over s+1}.
\label{eq:i3}
\ee
From \eqref{eq:i1}, \eqref{eq:i2} and \eqref{eq:i3}, we get
\beqas
&&\left\|{\partial\ue\over\partial t}(t)-{\partial u_0\over\partial t}(t)-\nabla_yw^r\left(\cdot,{\cdot\over\ep}\right)\left({\partial u_{0r}\over\partial t}(t)-g_{1r}\right)\right\|_H+\\
&&\qquad\qquad\|\curl\ue-\curl u_0-\scurl_yN^r\left(\cdot,{\cdot\over\ep}\right)\curl u_0(\cdot)_r\|_H\le c\ep^{s\over s+1}.
\eeqas
\eproof

\subsection{Corrector for multiscale problem}
For the case of more than two scales, we cannot deduce an explicit homogenization error. However, we can deduce correctors for the case where $\ep_{i-1}/\ep_i$ is an integer for all $i=2,\ldots,n$. We use the operator ${\cal T}_n^\ep$ and ${\cal U}_n^\ep$ defined as follows. 
We define the map
\[
{\cal T}_n^\ep(\phi)(x,\by)
=
\phi\Bigl(\ep_1\Bigl[{x\over\ep_1}\Bigr]
+
\ep_2\Bigl[{y_1\over\ep_2/\ep_1}\Bigr]
+
\ldots
+
\ep_n\Bigl[{y_{n-1}\over\ep_n/\ep_{n-1}}\Bigr]+\ep_ny_n\Bigr)
\]
for $\phi\in L^1(D)$ extended to 0 outside $D$.
Letting $D^{\ep_1}$ be the $2\ep_1$ neighbourhood of $D$, we have 
\be
\int_D\phi dx=\int_{D^{\ep_1}}\int_{Y_1}\cdots\int_{Y_n}{\cal T}_n^\ep(\phi)dy_n\cdots dy_1dx
\label{eq:Tepint}
\ee
for all $\phi \in L^1(D)$.
If a sequence $\{\phi^\ep\}_\ep$ $(n+1)$-scale converges to $\phi(x,y_1,\ldots,y_n)$, then 
\[
{\cal T}_n^\ep(\phi)\wc \phi(x,y_1,\ldots,y_n)
\]
in $L^2(D\times Y_1\times\ldots\times Y_n)$. We define the operator ${\cal U}_n^\ep$ 
\beqas
{\cal U}_n^\ep(\Phi)(x)
=
\int_{Y_1}\cdots\int_{Y_n}\Phi\Bigl(\ep_1\Bigl[{x\over\ep_1}\Bigr]
+
\ep_1 t_1,{\ep_2\over\ep_1}\Bigl[{\ep_1\over\ep_2}\Bigl\{{x\over\ep_1}\Bigr\}\Bigr]
+
{\ep_2\over\ep_1}t_2,\cdots,
\\
{\ep_n\over\ep_{n-1}}\Bigl[{\ep_{n-1}\over\ep_n}\Bigl\{{x\over\ep_{n-1}}\Bigr\}\Bigr]
+
{\ep_n\over\ep_{n-1}}t_n,\Bigl\{{x\over\ep_n}\Bigr\}\Bigr)dt_n\cdots dt_1
\eeqas
for all functions $\Phi\in L^1(D\times\bY)$. For each function $\Phi\in L^1(D\times\bY)$ we have 
\be
\int_{D^{\ep_1}}{\cal U}^\ep_n(\Phi)dx=\int_D\int_\bY\Phi(x,\by)d\by dx.
\label{eq:calUint}
\ee
The proofs for \eqref{eq:Tepint} and \eqref{eq:calUint} may be found in \cite{CDG}. We then have:
\be
{\cal T}^\ep_n\left({\partial\ue\over\partial t}\right)\wc {\partial u_0\over\partial t}+\sum_{i=1}^n{\partial\over\partial t}\nabla\fu_i,
\label{eq:wlTepduedt}
\ee
and
\be
{\cal T}^\ep_n(\curl\ue)\wc \curl u_0+\sum_{i=1}^n\scurl_{y_i}u_i
\label{eq:wlTepcurlue}
\ee
in $L^2(D\times\bY)$ when $\ep\to 0$. 

We have the following result.
\begin{theorem}\label{thm:mscorrector}
Assume that $g_0=0$, $g_1\in W$ and $f\in H^1(0,T;H)$. We have
\beqas
\lim_{\ep\to 0}\left\|{\partial\ue\over\partial t}-{\cal U}_n^\ep\left({\partial u_0\over\partial t}+\sum_{i=1}^n\nabla_{y_i}{\partial \fu_i\over\partial t}\right)\right\|_{L^\infty(0,T;H)}+
\left\|\curl\ue-{\cal U}^\ep_n\left(\curl u_0+\sum_{i=1}^n\scurl_{y_i} u_i\right)\right\|_{L^\infty(0,T;H)}=0.
\eeqas
\end{theorem}
\bproof
We consider
\beqas
E^\ep(t)&=&\int_D\int_\bY {\cal T}^\ep_n(b^\ep)\left({\cal T}^\ep_n\left({\partial\ue\over\partial t}\right)(t)-\left({\partial u_0\over\partial t}+\sum_{i=1}^n{\partial\over\partial t}\nabla\fu_i\right)(t)\right)\\
&&\qquad\qquad\qquad\qquad\cdot\left({\cal T}^\ep_n\left({\partial\ue\over\partial t}\right)(t)-\left({\partial u_0\over\partial t}+\sum_{i=1}^n{\partial\over\partial t}\nabla\fu_i\right)(t)\right)d\by dx\\ 
&&\qquad+\int_D\int_\bY{\cal T}^\ep_n(a^\ep)\left({\cal T}^\ep_n(\curl\ue)-\left(\curl_0+\sum_{i=1}^n\scurl_{y_i} u_i\right)\right)\\
&&\qquad\qquad\qquad\qquad\cdot \left({\cal T}^\ep_n(\curl\ue)-\left(\curl_0+\sum_{i=1}^n\scurl_{y_i} u_i\right)\right)d\by dx.
\eeqas
We note that
\beqas
&&\lim_{\ep\to 0}\int_D\int_\bY\left[{\cal T}^\ep_n(b^\ep){\cal T}^\ep_n\left({\partial\ue\over\partial t}\right)(t)\cdot{\cal T}^\ep_n\left({\partial\ue\over\partial t}\right)+{\cal T}^\ep_n(a^\ep){\cal T}^\ep_n(\curl\ue)\cdot{\cal T}^\ep_n(\curl\ue)\right]d\by dx\\
&&=\lim_{\ep\to 0}\int_D\left[b^\ep(x){\partial\ue\over\partial t}\cdot{\partial\ue\over\partial t}+a^\ep(x)\curl\ue(t,x)\cdot\curl\ue(t,x)\right]dx dt\\
&&=\lim_{\ep\to 0}\int b^\ep(x)g_1(x)\cdot g_1(x)dx+2\int_0^t\int_D f(t,x){\partial\ue\over\partial t}dxdt\\
&&=\int_D\left(\int_\bY b(x,\by)d\by\right)g_1(x)\cdot g_1(x)dx+2\int_0^t\int_D f(t,x){\partial u_0\over\partial t}dxdt
\eeqas
where we have used the energy formula for wave equation (see Lions and Magenes \cite{LionMagenes}) and the initial condition $g_0=0$. Using \eqref{eq:wlTepduedt} and \eqref{eq:wlTepcurlue}, we have 
\beqas
\lim_{\ep\to 0}E^\ep(t)&=&\int_D\left(\int_\bY b(x,\by)d\by\right)g_1(x)\cdot g_1(x)dx+2\int_0^t\int_D f(t,x){\partial u_0\over\partial t}(t,x)dxdt\\
&&-\int_D\int_\bY b(x,\by)\left({\partial u_0\over\partial t}+\sum_{i=1}^n{\partial\over\partial t}\nabla_{y_i}\fu_i\right)\cdot\left({\partial u_0\over\partial t}+\sum_{i=1}^n{\partial\over\partial t}\nabla_{y_i}\fu_i\right)d\by dx\\
&&-\int_D\int_\bY a(x,\by)(\curl u_0+\sum_{i=1}^n\scurl_{y_i}u_i)\cdot(\curl u_0+\sum_{i=1}^n\scurl_{y_i}u_i).
\eeqas 
From \eqref{eq:Gn} and \eqref{eq:tildeGn}, we have
\beqas
&&\int_D\int_\bY b(x,\by)\left({\partial u_0\over\partial t}+\sum_{i=1}^n{\partial\over\partial t}\nabla_{y_i}\fu_i\right)\cdot\left({\partial u_0\over\partial t}+\sum_{i=1}^n{\partial\over\partial t}\nabla_{y_i}\fu_i\right)d\by dx\\
&&=\int_D\int_\bY b(x,\by)\left[{\partial\over\partial t}(u_0+\sum_{i=1}^{n-1}\nabla_{y_i}\fu_i)+\left({\partial\over\partial t}(u_0+\sum_{i=1}^{n-1}\nabla_{y_i}\fu_i)_k-g_{1k}\right)\nabla_{y_n}w^k_n\right]\\
&&\qquad\qquad\cdot
\left[{\partial\over\partial t}(u_0+\sum_{i=1}^{n-1}\nabla_{y_i}\fu_i)+\left({\partial\over\partial t}(u_0+\sum_{i=1}^{n-1}\nabla_{y_i}\fu_i)_k-g_{1k}\right)\nabla_{y_n}w^k_n\right]\\
&&=\int_D\int_\bY b(x,\by)\left[{\partial\over\partial t}(u_0+\sum_{i=1}^{n-1}\nabla_{y_i}\fu_i)_k(e^k+\nabla_{y_n}w^k_n)-g_{1k}\nabla_{y_n}w^k_n\right]\\
&&\qquad\qquad\cdot \left[{\partial\over\partial t}(u_0+\sum_{i=1}^{n-1}\nabla_{y_i}\fu_i)_l(e^l+\nabla_{y_n}w^l_n)-g_{1l}\nabla_{y_n}w^l_n\right].
\eeqas
From \eqref{eq:cellbn} and \eqref{eq:bnminus1}, this equals
\[
\int_D\int_{\bY_{n-1}}b^{(n-1)}_{kl}{\partial\over\partial t}(u_0+\sum_{i=1}^{n-1}\nabla_{y_i}\fu_i)_k{\partial\over\partial t}(u_0+\sum_{i=1}^{n-1}\nabla_{y_i}\fu_i)_ld\by_{n-1}dx+\int_D\int_\bY b(x,\by)\nabla_{y_n}w^k_n\cdot\nabla_{y_n}w^l_n g_{1k}g_{1l}.
\]
Continuing this process, this expression equals
\[
\int_D b^0(x){\partial u_0\over\partial t}\cdot{\partial u_0\over\partial t}dx+\sum_{i=1}^n\int_D\int_{\bY_i} b_i(x,\by_i)\nabla_{y_i}w^k_i\cdot\nabla_{y_i}w^l_ig_{1k}g_{1l}d\by_i dx.
\]
On the other hand we have
\[
\int_D\int_\bY a(x,\by)(\curl u_0+\sum_{i=1}^n\scurl_{y_i}u_i)\cdot(\curl u_0+\sum_{i=1}^n\scurl_{y_i}u_i)d\by dx=\int_D a^0(x)\curl u_0\cdot\curl u_0dx.
\]
Therefore
\beqas
\lim_{\ep\to 0}E^\ep(t)=\int_D\left(\int_\bY b(x,\by)d\by\right)g_1\cdot g_1 dx+2\int_0^t\int_Df{\partial u_0\over\partial t}dxdt-\int_D b^0(x){\partial u_0\over\partial t}\cdot{\partial u_0\over\partial t}dx\\
-\sum_{i=1}^n\int_D\int_\bY b_i(x,\by_i)\nabla_{y_i}w^k_i\cdot\nabla_{y_i}w^l_ig_{1k}g_{1l}d\by_i dx-\int_D a^0(x)\curl u_0\cdot\curl u_0dx.
\eeqas
From \eqref{eq:homogenizedeq}, we get
\[
\int_D b^0(x){\partial u_0\over\partial t}\cdot{\partial u_0\over\partial t}dx+\int_D a^0(x)\curl u_0\cdot\curl u_0dx=\int b^0(x)g_1\cdot g_1 dx+2\int_0^t\int_D f{\partial u_0\over\partial t}dxdt.
\]
From \eqref{eq:cellbi} we have
\[
\int_{Y_i}b^i(x,\by_i)\nabla_{y_i}w^k_i\cdot\nabla_{y_i}w^l_i dy_i=-\int_{Y_i}b^i(x,\by_i)e^k\cdot\nabla_{y_i}w^l_i.
\]
Thus
\beqas
&&\int_D\left(\int_\bY b(x,\by)d\by\right)g_1\cdot g_1 dx-\sum_{i=1}^n\int_D\int_{\bY_i} b^i(x,\by_i)\nabla_{y_i}w^k_i\cdot\nabla_{y_i}w^l_ig_{1k}g_{1l}d\by_i dx\\
&&\qquad\qquad=\int_D\left(\int_\bY b(x,\by)d\by\right)g_1\cdot g_1 dx+\sum_{i=1}^n\int_D\int_{\bY_i}b^i(x,\by_i)e^k\cdot\nabla_{y_i}w^l_ig_{1l}g_{1k}d\by_i dx.
\eeqas
We consider
\beqas
&&\int_D\left(\int_\bY b(x,\by)d\by\right)g_1\cdot g_1 dx+\int_D\int_\bY b(x,\by)e^k\cdot\nabla_{y_n}w^l_n g_{1l}g_{1k} d\by dx\\
&&=\int_D\left(\int_\bY b_{ij}(x,\by)d\by\right)g_{1j}(x)\cdot g_{1i}(x) dx+\int_D\int_\bY b_{si}{\partial w^j_n\over\partial y_s}(x,\by)g_{1j}(x)g_{1i}(x)dx\\
&&=\int_D\int_{\bY_{n-1}}\left(\int_{Y_n}b_{is}(\delta_{js}+{\partial w^j_n\over\partial y_s})dy_n\right)g_{1j}(x)g_{1i}(x) d\by_{n-1}dx\\
&&=\int_D\int_{\bY_{n-1}}b^{n-1}_{ij}(x,\by_{n-1})g_{1j}(x)g_{1i}(x)dxd\by_{n-1}.
\eeqas
Continuing this, we get 
\[
\int_D\left(\int_\bY b(x,\by)d\by\right)g_1\cdot g_1 dx+\sum_{i=1}^n\int_D\int_{\bY_i}b^i(x,\by_i)e^k\cdot\nabla_{y_i}w^l_ng_{1l}g_{1k}d\by_i dx=\int_D b^0(x)g_1(x)\cdot g_1(x) dx.
\]
Thus
\beqas
\lim_{\ep\to 0}E^\ep(t)=0.
\eeqas
We show that the convergence is uniform. To make the notation concise, we denote by
\[
A(t)={\cal T}^\ep_n\left({\partial\ue\over\partial t}\right)(t)-\left({\partial u_0\over\partial t}+\sum_{i=1}^n{\partial\over\partial t}\nabla\fu_i\right)(t).
\]
We consider
\beqas
\left|\int_D\int_\bY \left[{\cal T}^\ep_n(b^\ep)A(t)\cdot A(t)-{\cal T}^\ep_n(b^\ep) A(s)\cdot A(s)\right]d\by dx\right|=\left|\int_D\int_\bY {\cal T}^\ep_n(b^\ep)(A(t)-A(s))\cdot(A(t)+A(s))d\by dx\right|\\
\le c\|A(t)-A(s)\|_{L^2(D\times\bY)^3}\|A(t)+A(s)\|_{L^2(D\times\bY)^3}.
\eeqas
We note that $\|{\partial\ue\over\partial t}\|_{L^\infty(0,T;H)}$ is uniformly bounded for all $\ep$, ${\partial u_0\over\partial t}\in L^\infty(0,T;H)$, and from \eqref{eq:dtfui}, ${\partial\over\partial t}\nabla_{y_i}\fu_i\in L^\infty(0,T;L^2(D\times\bY_i)^3)$. Further 
\beqas
A(t)-A(s)=\int_s^t\left[{\cal T}^\ep_n\left({\partial^2\ue\over\partial t^2}\right)(\tau)-\left({\partial^2u_0\over\partial t^2}(\tau)+\sum_{i=1}^n{\partial^2\over\partial t^2}\nabla_{y_i}\fu_i(\tau)\right)\right]d\tau
\eeqas
so
\begin{eqnarray}
&&\|A(t)-A(s)\|_{L^2(D\times\bY)^3}\nonumber\\
&&\qquad\qquad\le \int_s^t\left[\left\|{\cal T}^\ep_n\left({\partial^2\ue\over\partial t^2}\right)(\tau)\right\|_{L^2(D\times\bY)^3}+\left\|{\partial^2u_0\over\partial t^2}(\tau)\right\|_{L^2(D)^3}+\sum_{i=1}^n\left\|{\partial^2\over\partial t^2}\nabla_{y_i}\fu_i(\tau)\right\|_{L^2(D\times\bY)^3}\right]d\tau
\label{eq:X2}
\end{eqnarray}
From \eqref{eq:duepdt}, we have that ${\partial^2\ue\over\partial t^2}$ is uniformly bounded in $L^2(0,T;H)$. By a similar argument using the compatible initial condition, we show that ${\partial^2u_0\over\partial t^2}\in L^2(0,T;H)$ which implies that ${\partial^2\over\partial t^2}\nabla_{y_i}\fu_i\in L^2(0,T;L^2(D\times\bY_i)^3)$. We then have
\beqas
\int_s^t\left\|{\cal T}^\ep_n\left({\partial^2\ue\over\partial t^2}\right)(\tau)\right\|_{L^2(D\times\bY)^3}d\tau\le (t-s)^{1/2}\left(\int_0^t\left\|{\cal T}^\ep_n\left({\partial^2\ue\over\partial t^2}\right)(\tau)\right\|^2_{L^2(D\times\bY)^3}\right)^{1/2}\\
\le c(t-s)^{1/2}\left(\int_0^t\left\|{\partial^2\ue\over\partial t^2}(\tau)\right\|_{L^2(D)^3}^2\right)^{1/2}\le c(t-s)^{1/2}.
\eeqas
By the same argument, we have similar estimates for other terms in \eqref{eq:X2}. We can perform similarly for the other terms in $E^\ep(t)$. From the Arzel\`a-Ascoli theorem, we deduce that $E^\ep(t)$ converges to 0 when $\ep\to 0$ uniformly for all $t\in [0,T]$. The conclusion of the proposition follows from the fact that
\beqas
\left\|{\partial\ue\over\partial t}(t)-{\cal U}_n^\ep\left({\partial u_0\over\partial t}+\sum_{i=1}^n{\partial\over\partial t}\nabla_{y_i}\fu_i\right)(t)\right\|_{L^2(D)^3}\le \left\|{\cal T}^\ep_n\left({\partial\ue\over\partial t}\right)(t)-\left({\partial u_0\over\partial t}+\sum_{i=1}^n{\partial\over\partial t}\nabla_{y_i}\fu_i\right)(t)\right\|_{L^2(D\times\bY)^3},
\eeqas
and 
\beqas
\left\|\curl\ue(t)-{\cal U}^\ep_n\left(\curl u_0+\sum_{i=1}^n\scurl_{y_i}u_i\right)(t)\right\|_{L^2(D)^3}\le \left\|{\cal T}^\ep_n(\curl\ue(t))-\left(\curl u_0+\sum_{i=1}^n\scurl_{y_i}u_i\right)(t)\right\|_{L^2(D)^3}.
\eeqas
\eproof


{\bf Acknowledgement} The authors gratefully acknowledge a postgraduate scholarship of Nanyang Technological University, the  AcRF Tier 1 grant RG30/16, the Singapore A*Star SERC grant 122-PSF-0007 and the AcRF Tier 2 grant MOE 2013-T2-1-095 ARC 44/13. 

\bibliographystyle{plain}
\bibliography{references}
\end{document}